\documentclass[11pt,english]{article}
\usepackage[T1]{fontenc}
\usepackage[latin9]{inputenc}
\usepackage{geometry}
\usepackage{babel}

\usepackage{amsmath}
\usepackage{amsfonts}
\usepackage{amssymb}
\usepackage{graphicx}
\usepackage{color, graphics}
\usepackage{natbib}
\usepackage{relsize}   
\usepackage{bbm}

\newtheorem{thm}{Theorem}[section]
\newtheorem{lemma}[thm]{Lemma}
\newtheorem{prop}[thm]{Proposition}
\newtheorem{cor}[thm]{Corollary}
\newtheorem{defn}[thm]{Definition}

\newtheorem{exa}[thm]{Example}

\newcommand{\qed}{\nobreak \ifvmode \relax \else
      \ifdim\lastskip<1.5em \hskip-\lastskip
      \hskip1.5em plus0em minus0.5em \fi \nobreak
      \vrule height0.75em width0.5em depth0.25em\fi}

\newcommand{\eqdist}{\overset{d}{=}}
\newcommand{\calG}{\mathcal{G}}
\newcommand{\calF}{\mathcal{F}}
\newcommand{\iidsim}{\overset{iid}{\sim}} 
\newcommand{\indsim}{\overset{ind}{\sim}} 
\newcommand{\idsim}{\overset{id}{\sim}} 

\newcommand{\Polya}{P\'{o}lya }
\newcommand{\Norm}{\text{N}} 
\newcommand{\E}{\mathbb{E}} 

\newcommand{\bP}{{\mathbb{P}}}
\newcommand{\bX}{{\mathbf{X}}}

\newcommand{\Times}{\mathlarger{\mathlarger{\mathlarger{\times}}}}

\begin{document}

\title{On a notion of partially conditionally identically distributed  sequences}

\author{Sandra Fortini\footnote{Corresponding author. Department of Decision Sciences, Bocconi University, Via Roentgen 1, 	20123 Milano, Italy. E-mail address: sandra.fortini@unibocconi.it
	}, Sonia Petrone, and Polina Sporysheva
       \\
       Bocconi University, Milan, Italy
       }

\maketitle


\begin{abstract}
A notion of conditionally identically distributed (c.i.d.) sequences has been studied as a form of stochastic dependence weaker than exchangeability, but  equivalent to it in the presence of stationarity.  
We extend such notion to families of sequences. Paralleling the extension from exchangeability to partial exchangeability in the sense of de Finetti, we propose a notion of {\em partially c.i.d.} dependence, which is shown to be equivalent to partial exchangeability for stationary processes. Partially c.i.d. families of sequences preserve attractive limit properties of partial exchangeability, and are asymptotically partially exchangeable. Moreover, we provide strong laws of large numbers and two central limit theorems. Our focus is on the asymptotic agreement of predictions and empirical means, which lies at the foundations of Bayesian statistics.
Natural examples of partially c.i.d. constructions are  interacting randomly reinforced processes satisfying certain conditions on the  reinforcement. 

\end{abstract}


\section{Introduction.} \label{sec:intro}

Exchangeability is a central notion in many areas of probability and related fields; we  refer to \cite{Kingman78exchangeability}, \cite{Aldous85}, \cite{austin2008}, \cite{kallenberg2005},  \cite{aldous2010} for classical, wide references. In Bayesian statistics, exchangeability is the fundamental probabilistic structure at the basis of learning, expressing the subjective probabilistic description of repeated experiments under similar conditions.
Exchangeable sequences are conditionally independent and identically distributed (i.i.d.). 

However, forms of competition, selection, and other sources of non stationarity, may break exchangeability, although the system may converge, asymptotically, to an exchangeable steady state. Thus, weaker notions of stochastic dependence, that do not assume  stationarity, yet preserve some main asymptotic properties of exchangeable processes, are of theoretical and applied interest. Based on results by \cite{Kallenberg88}, \cite{BertiPratelliRigo2004LimThForIID} introduce a notion of {\em conditionally identically distributed} (c.i.d.) sequences, as a form of stochastic dependence weaker than exchangeability but equivalent to it  for stationary sequences.  Roughly speaking, a sequence of random variables $(X_n)_{n \geq 1}$ is  c.i.d. if, for every $n \geq 0$, $X_{n+1},X_{n+2},\dots$ are conditionally identically distributed, given the past $X_1, \ldots, X_n$
(for $n=0$ the property reads $X_1,X_2,\dots$ are identically distributed). A precise definition is given in Section \ref{sec:cid}. These processes are the starting point of our study.

Notions of partial exchangeability are needed for more complex phenomena, which  can still be described by some form of probabilistic invariance, under specific subclasses of permutations. See \cite{diaconisFreedman1980}, \cite{kallenberg2005},  \cite{aldous2010}. 
A basic notion is partial exchangeability in the sense of de Finetti \citep{deFinetti1937-scparziale} (called internal exchangeability for a family of sequences by \cite{Aldous85}). 
In  Bayesian statistics, partial exchangeability in the sense of de Finetti is the fundamental probabilistic dependence behind inference for multiple experiments. Roughly speaking, observations are exchangeable within each experiment, but not across experiments; the probabilistic dependence among  sequences allows borrowing strength across experiments. Throughout this paper, by partial exchangeability we mean partial exchangeability in the sense of de Finetti.
Again, different forms of non-stationarity may break the symmetry of partial exchangeability. It seems natural to ask how the notion of c.i.d. sequences can be extended to a notion of {\em partially} c.i.d. processes, the way that partial exchangeability extends exchangeability.
Such extension is the main objective of the present work.

We introduce a notion of {\em partially c.i.d.} families of sequences that is shown to be equivalent to partial exchangeability under stationarity. Then, we prove that partially c.i.d. sequences preserve some main limit properties of partially exchangeable sequences. In particular, the joint predictive distributions and the joint empirical distributions converge (weakly) to the same random limit, almost surely. Moreover, partially c.i.d. sequences are asymptotically partially exchangeable.
The asymptotic agreement of frequencies and predictions is of fundamental interest in Bayesian statistics, where probability has a subjective interpretation, showing the frequentist basis of the subject's probabilistic learning. Such agreement is ensured (in the subject's opinion) for exchangeable and  partially exchangeable sequences. Our result shows that it is still ensured when relaxing the assumption of stationarity from partial exchangeability.
Marginally, these results are not surprising, as partially c.i.d. sequences are marginally c.i.d. and the limit behavior of c.i.d. sequences has been studied (\cite{BertiPratelliRigo2004LimThForIID}, \cite{BertiPratelliRigo2012LimThforEmpirProcBasedOnDepData}). Yet, for multiple sequences, the {\em joint} limit behavior is not obvious, as the sequences are stochastically dependent. Notice that they remain asymptotically dependent, if the random marginal limit measures are dependent.

These limit results are refined in Section \ref{sec:SLLN}, where we provide a strong law of large numbers for partially c.i.d. sequences, and in Section \ref{sec:CLT}, where we give two central limit theorems, 
for the scaled cumulative forecast errors and the scaled difference between empirical means and predictions, respectively.
Beyond fundamental issues, the possibility of approximating  predictions by  empirical means, with an approximation error given by a central limit theorem,
can be of interest for hypothesis testing and model checking as well as for facilitating computations in Bayesian prediction with large sample size.
\medskip

Areas of applications include interacting evolutionary phenomena that, while not being stationary, tend towards an equilibrium state of partial exchangeability. We provide several examples throughout the paper. In particular, a natural class of partially c.i.d. families of sequences are interacting randomly reinforced processes satisfying certain conditions on the random reinforcement (see Section \ref{sec:partialExch}). 
In the case of a single sequence, contructions based on a time-varying or random reinforcement (see \cite{pemantle2007}) are indeed main examples of c.i.d. sequences. The time-dependent urn scheme of \cite{pemantle1990time} is a basic case, and further examples are Randomly Reinforced Urns (RRUs), a special case of Generalized \Polya Urns (see  \cite{athreyaNey1972}, and \cite{pemantle2007}), with a diagonal and random replacement matrix.
In the classical two-color \Polya urn, a ball is drawn at each step, and returned to the urn together with an additional ball of the same color. Then 
the sequence of colors extracted at successive draws is exchangeable.
However, the number of additional balls placed in the urn at each step, the \emph{reinforcement}, may be random, making the sequence of colors no longer exchangeable. In fact, if the random reinforcements are independent of the color extracted, the generated sequence of colors is c.i.d. Areas of application include adaptive clinical trails (\cite{huRosenbeger2006}, \cite{antogniniGiovagnoli2015}), sequential design, two armed bandit problems and reinforcement learning (\cite{beggs2005}) and Bayesian inference.  The two-color randomly reinforced urn scheme proposed by \cite{durhamYu1990}
for response-adaptive clinical trials generates a c.i.d. sequence under the null hypothesis of equivalence of the two treatments. Extensions and theoretical properties are given, among others, by
\cite{durhamFlournoyLi1998}, \cite{Muliere2006RRU}, \cite{mayFlournoy2009}; see \cite{flournoyMaySecchi2012} for a review.
The notion of partially c.i.d. sequences extends the spectrum of application to multiple experiments, for example clinical trials in multiple centers. There is an increasing interest in interacting randomly reinforced processes in many fields (see e.g. \cite{paganoniSecchi2004};  \cite{crimaldi2015fluctuation}).
Yet,  the literature is somehow fragmented, as many results are tailored for specific constructions and aims; the notion of partially c.i.d. processes may provide a useful reference framework.

Processes with reinforcement are the basis of many important constructions in Bayesian nonparametric inference; extensions to a random reinforcement are of interest also in this area (\cite{Bassetti2010CidSSS}, \cite{airoldiCosta2014}).
Applications to competitive networks are shown by \cite{caldarelli2013} and c.i.d. sequences arise as a particular case of the generalized Indian Buffet process \citep{bertiIBP2015} for competitive feature selection.
Multivariate extensions of these constructions do not appear to have been studied. In fact, a problem of interest in Bayesian nonparametrics is to characterize partially exchangeable families of reinforced processes with a tractable form of the predictive laws. However, this problem is somehow unsolved (see the discussion 
in  \cite{wadeMongelluzzoP2011}, \cite{Lee2013} and \cite{bacalladoDiaconisHolmes2015}). Natural constructions of dependent reinforced processes fail to be partially exchangeable. However, they may generate partially c.i.d.,  therefore asymptotically partially exchangeable, sequences. 
 Implicitly, these partially c.i.d. constructions characterize novel classes of prior distributions for dependent random measures, as the weak limits of the joint predictive distributions. However, finding explicit expressions of these prior distributions is difficult. 
Some results for c.i.d. or \lq \lq quasi-c.i.d.'' sequences are in \cite{durhamFlournoyLi1998} and \cite{alettiMaySecchi2007,{alettiMaySecchi2012}, {alettiGhigliettiPaganoni2013}}. 
 \cite{pemantle1990} gives non-convergence results that apply to partially c.i.d. sequences. The potential of stochastic approximation methods for randomized urn schemes is further highlighted by \cite{laruelle2013randomized}. A study of the explicit limit law for partially c.i.d. sequences is beyond the aim of this paper.  Yet, we show an example, with dependent Gaussian sequences, for which the limit law can be computed and has a parametric expression (Section \ref{sec:partialExch}).

We start our study with a brief review of c.i.d. sequences, including  some facts that appear to be novel, in Section \ref{sec:cid}. 
We introduce the notion of partially c.i.d. sequences in Section \ref{sec:pcid} and prove some main limit results in Section  \ref{sec:partialExch}.
A strong law of large numbers  and two central limits theorems are presented in Sections \ref{sec:SLLN} and \ref{sec:CLT}. Several examples are provided throughout the paper. 

\section{Exchangeability, spreadability and c.i.d. sequences.}
\label{sec:cid}

Let  $(X_n)_{n \geqslant 1}$  be a sequence of random variables defined on a probability space  $(\Omega,\mathcal{F},\mathbb{P})$, with each $X_{n}$ taking values in a Polish space $\mathbb{X}$, endowed with the Borel $\sigma$-algebra  $\mathcal{X}$. 
Spaces of probability measures are endowed with the topology of weak convergence and with the sigma-algebra generated by the evaluation maps. A random probability measure is a random element taking values in a space of probability measures.
Conditional distributions refer to 
regular versions; equalities of random probability measures hold $\bP$-almost surely ($\bP$-a.s.). Additionally, we use $X_{1:n}$ to denote the vector $(X_1, \ldots, X_n)$; $(X_n)$ to denote the sequence  $(X_n)_{n \geqslant 1}$ and $\eqdist$ to denote equality in distribution.

The sequence $(X_n)$ is exchangeable if its probability law is invariant under any finite permutation, that is, if
$$ (X_1, X_2, \ldots) \eqdist (X_{\pi(1)}, X_{\pi(2)}, \ldots) 
$$
for any permutation $\pi$ of $\mathbb N$ for which $\pi(n)=n$, except for a finite number of $n$'s. 

An equivalent notion of invariance under subsequence selection, or spreading invariance, is discussed in \cite{Kingman78exchangeability} and  \cite{Aldous85}. Following the terminology of \cite{Kallenberg88}, we say that a sequence $(X_n)$ is {\em spreadable} if
$$ (X_1, X_2, \ldots) \eqdist (X_{k_1}, X_{k_2}, \ldots) \quad \mbox{for every } k_1 < k_2 < \cdots.
$$
It is easy to show that an exchangeable sequence is spreadable, and the converse implication is proved by \cite{ryllNad1957}. By de Finetti's  representation theorem, the law of an exchangeable sequence is a mixture of distributions of i.i.d. random variables (\cite{Aldous85}, Part I, Section 3). 
\cite{ryllNad1957} proves that the same conclusion holds for spreadable sequences. In fact, \cite{Kallenberg88} notices that the  representation theorem is a consequence of the mean ergodic theorem,
and can be proved for stationary sequences that satisfy a condition weaker than spreadability. 
\begin{prop} \label{prop:kallenberg}
{\em \cite[Proposition 2.1]{Kallenberg88}}.
A stationary sequence $(X_n)$  that  satisfies 
\begin{equation} \label{eq:mgspreadability}
(X_1, \ldots, X_n, X_{n+1}) \eqdist (X_1, \ldots, X_n, X_{n+k}) \quad \mbox{for all integers $k \geq 1$ and $n \geq 1$} 
\end{equation}
is exchangeable.
\end{prop}
Notice that (\ref{eq:mgspreadability}) is equivalent to  $(X_1, \ldots, X_n, X_{n+1}) \eqdist (X_1, \ldots, X_n, X_{n+2})$ for all $ n \geq 1$. Condition (\ref{eq:mgspreadability}) says that all future observations are conditionally identically distributed, given the past.
 Extending this notion, \cite{BertiPratelliRigo2004LimThForIID}  give the following definition. Assume that the sequence $(X_n)$ is adapted to a filtration $\calG =\left(\mathcal{G}_n\right)_{n\geqslant 0}$. Then,  $(X_{n})$ is \textit{conditionally identically distributed with respect to the filtration $\mathcal{G}$}, or $\mathcal{G}$-c.i.d. whenever
\begin{equation}
\begin{aligned} \label{eq:naturalCid-rev}
&\mathbb E[ f(X_{n+k}) \mid \mathcal G_n] = \mathbb E[ f(X_{n+1}) \mid \mathcal G_n]\quad \mathbb P\mbox{-a.s.} \\
&\mbox{for all $k \geq 1$, $n \geq 0$ and all bounded measurable functions $f: \mathbb{X} \rightarrow \mathbb{R}$} \nonumber. 
\end{aligned}
\end{equation}
Roughly speaking, the future observations $X_{n+k}$ are identically distributed, given the past $\calG_n$. 
In particular, the  $X_n$ are marginally identically distributed.
 When considering the natural filtration, i.e.  when $\calG_n$ is the
sigma-field  generated by $(X_1, \ldots, X_n)$ and $\mathcal G_0=\{\emptyset,\Omega\}$, $(X_n)$ is called {\em naturally c.i.d.}, or, simply, c.i.d.,  omitting the filtration. If $(X_n)$ is $\cal G$-c.i.d., then $(X_n)$ is also c.i.d. with respect to any coarser filtration to which it is adapted; in particular, it is c.i.d.  An exchangeable sequence is clearly c.i.d.

The following proposition (\cite{Kallenberg88} Proposition 2.2,
and \cite{BertiPratelliRigo2004LimThForIID}) gives equivalent conditions.
\begin{prop} \label{prop:kallenberg2}
The following properties (i), (ii) and (iii) are equivalent:
\begin{itemize}
\item[(i)] The sequence $(X_n)$ is $\calG$-c.i.d.
\item[(ii)] 
The process of  predictive measures $Q_n(\cdot) := \mathbb{P}[X_{n+1} \in \cdot \mid \calG_{n}]$ is a measure-valued $\mathcal{G}$-martingale.
\item[(iii)]
For each finite $\calG$-stopping time $\tau$, $X_{\tau+1} \eqdist X_1$.
\end{itemize}
\end{prop}
The martingale condition (ii) asserts that  $(Q_n)$ is a sequence of random probability measures satisfying
	$\mathbb E[Q_{n+1}(A)\mid\mathcal G_n]=Q_n(A)$, $\mathbb P\mbox{-a.s.} $, for every $n\geq 0$ and every $A\in \mathcal X$
	(\cite{horowitz1985}); this  is equivalent to:
	\begin{equation} \label{eq:CID-II}
\left(\mathbb{E}\left[f(X_{n+1})\mid\mathcal{G}_{n}\right]\right)_{n\geqslant 0}\;\text{is a $\mathcal{G}$-martingale}, 
 \text{ for all bounded measurable functions $f: \mathbb{X} \rightarrow \mathbb{R}$.}
\end{equation}
 
We present below some facts, that appear to be novel, which provide further insights into the connections with the notion of exchangeability.
Further characterizations are given by  \cite{BertiPratelliRigo2004LimThForIID,  BertiPratelliRigo2012LimThforEmpirProcBasedOnDepData}.   

An exchangeable sequence can be characterized by the sequence of predictive distributions $(Q_n)_{n\geq 0}$, with $Q_0:=\mathbb P[X_1\in\cdot\;]$ and $Q_n:=\mathbb P[X_{n+1}\in\cdot\mid X_{1:n}]$, for $n\geq 1$ \citep{fortiniLadelliRegazzini2000}. In particular, a necessary condition for $(Q_n)$ to define an exchangeable probability law for $(X_n)$ is that $Q_{n}$ is a symmetric function of $(X_1, \ldots, X_n)$. For c.i.d. sequences, the symmetric role of past observations in prediction is lost. Indeed, a c.i.d. sequence $(X_n)$ is exchangeable if and only if $\mathbb{P}\left[X_{n+1}\in \cdot\mid X_1\in A_1,\dots,X_n\in A_n\right]$ is  symmetric in $(A_1,\dots,A_n)$, for every $n\geq 1$ and every $A_1,\dots, A_n$ with $\mathbb P[X_1 \in A_1,\dots, X_n \in A_n]>0$.
This result is proved as a Corollary of Proposition \ref{prop:predCharacterization} in Section \ref{sec:pcid}.

\medskip

A second fact is related to the lack of a de Finetti-type representation theorem for c.i.d. sequences. If $(X_n)$ is exchangeable, by de Finetti's representation theorem there exists a random probability measure $\alpha$ such that, conditionally on $\alpha$,  $X_1,X_2,\dots$ are i.i.d, with distribution $\alpha$;  $\alpha$ is called the \textit{directing random measure} of $(X_n)$ and coincides $\mathbb P$-a.s. with the weak limit of both the sequence of the empirical distributions and the sequence of the predictive distributions of $(X_n)$ (\cite{Aldous85}).
Thus, although an exchangeable sequence $(X_n)$ may describe an {\em evolutionary}  process,  
the representation theorem implies that $(X_n)$ is probabilistically equivalent to {\em static} random sampling from its directing measure.
The possible lack of stationarity clearly implies that no equivalence with a static phenomenon, that is, no similar representation result, is possible for c.i.d. sequences. Intuitively, $X_n$ is sampled from an evolving population that has unpredictable dynamics and converges to a random steady state. The following proposition gives a state-space like construction.

\begin{prop}
	Consider a process $(X_n, F_n)$ where $F_n$ are 
	random probability measures on $\mathbb X$. If 
the following conditions hold:
	\begin{itemize}
			\item[(i)] conditionally on $(F_n)$, the $X_n$ are independent and the conditional distribution of $X_n$ is $F_n$;
	\item[(ii)] $(F_n)$ is a measure-valued martingale with respect to its natural filtration, 
\end{itemize}
then the process $(X_n)$ is c.i.d.
	\end{prop}

\noindent {\sc Proof}. 
We prove that, under (i) and (ii), $(X_n)$ is $\calG$-c.i.d., where $\mathcal G_0=\{\emptyset,\Omega\}$ and $\calG_n$ is the sigma-field generated by $(X_{1:n}, F_{1:n})$. Therefore,  it is also c.i.d. 
By (i) and (ii), $\mathbb P[X_k\in\cdot\;]=\mathbb E[F_k(\cdot)]=\mathbb E[F_1(\cdot)]=\mathbb P[X_1\in\cdot\;]$. Furthermore, by
(i), for every $n\geq 1$, every bounded continuous functions $g,g_1,\dots, g_n$ and every $A\in\mathcal X$,
$$
\mathbb E[g_1(X_1)g_2(X_2)\dots,g_n(X_n)g(F_{n+k}(A))\mid F_{1:n}]=   \mathbb E[g(F_{n+k}(A))\mid F_{1:n}] \; \prod_{s=1}^n\mathbb E[g_s(X_s)\mid F_s] \; ; 
$$
therefore, $F_{n+k}(A)$ is conditionally 
independent of $X_{1:n}$, given $F_{1:n}$. Then, for any $k \geq 1$, $ n\geq 1$ and $A\in\mathcal X$,
\begin{eqnarray*} 
\mathbb P[X_{n+k} \in A \mid X_{1:n}, F_{1:n}] &=&
\mathbb E [ \mathbb P[X_{n+k} \in A \mid X_{1:n}, (F_n)] \mid X_{1:n}, F_{1:n}]
=\mathbb E[F_{n+k}(A) \mid  F_{1:n}] \\
&=& \mathbb E[F_{n+1}(A) \mid  F_{1:n}] 
=\mathbb E[F_{n+1}(A) \mid X_{1:n}, F_{1:n}]
= \mathbb P[X_{n+1} \in A \mid X_{1:n}, F_{1:n}]
\end{eqnarray*}
where the third equality follows from the martingale property (ii). Thus, $(X_n)$ is $\calG$-c.i.d.\\
\noindent $\square$  


Abusing notation, we will write $X_n \mid F_n \indsim F_n$ to denote the dependence structure under (i) above. The following example, which elaborates from Example 1.3 in \cite{BertiPratelliRigo2004LimThForIID}, suggests that a c.i.d. process may have a state-space representation in terms of a finite-dimensional latent process $(\theta_n)$ converging to a random limit, a result that we do not pursue further here.

\begin{exa} \label{ex:gaussian cid}
{\em Let $(X_n, \theta_n)$ be described by the following equations:
\begin{eqnarray} \label{eq:state-space}
X_n &=& \theta_n + \epsilon_n , \quad \epsilon_n \indsim \Norm(0, c-b_n); \nonumber \\
\theta_n &=& \theta_{n-1} + v_n , \quad v_n \indsim \Norm(0, b_n-b_{n-1}),\quad \quad  n\geq 1,
\end{eqnarray}
for some $0=b_0 < b_1 < b_2 < \cdots < c$, with $b_n \rightarrow c' < c$, and with  $\theta_0$, $(\epsilon_n)$ and $(v_n)$ independent. For brevity, let  $\theta_0$ be a fixed constant. In other words,
$$ X_n \mid \theta_n \indsim \Norm(\theta_n, c-b_n) $$
and $(\theta_n)$ is the damped random walk described by the state equation  (\ref{eq:state-space}). It is easy to show that the process $(X_n)$ is c.i.d.,  with $X_n \idsim \Norm(\theta_0, c)$ (where $\idsim$ means identically distributed). 
Some computations show that $\theta_n=\theta_0+\sum_{i=1}^n v_i$ converges in distribution to a random limit $\theta \sim \Norm(\theta_0, c') $ 
(indeed, $(\theta_n)$ is a uniformly integrable martingale, therefore it converges a.s. to a random limit $\theta$). Moreover,  the predictive distribution of $X_{n+1}$, given $X_{1:n}$, converges to a $\Norm(\theta, c-c')$. By Lemma 8.2 in \cite{Aldous85}, this implies that $(X_n)$ is asymptotically exchangeable, and the exchangeable limit law has directing random measure given by $\Norm(\theta, c-c')$; that is, roughly speaking,  $X_n \mid \theta \iidsim \Norm(\theta, c-c')$ for large $n$.
}
\end{exa}

Although generally not exchangeable, $\mathcal{G}$-c.i.d. sequences preserve some  attractive limit properties of exchangeable sequences \citep {BertiPratelliRigo2004LimThForIID}. In particular, $\mathbb{P}$-a.s., the empirical distributions $\sum_{i=1}^n\delta_{X_i}(\cdot)/n$ and the predictive distributions $\mathbb{P}[X_{n+1}\in \cdot \mid\mathcal{G}_n]$ converge  to the same random probability measure $\alpha$.
Moreover, $\mathbb{P}[X_{n+1}\in \cdot \mid\mathcal{G}_n]=\mathbb{E}[\alpha(\cdot)\mid\mathcal{G}_n]$.
Convergence of the predictive distributions to $\alpha$ implies that the sequence $(X_n)$ is asymptotically exchangeable and the exchangeable limit law has directing measure $\alpha$; that is, $\left(X_{n+1},X_{n+2},\dots\right)\overset{d}{\rightarrow}\left(Z_1,Z_2,\dots\right)$ as $n\rightarrow\infty$, for some exchangeable sequence $(Z_{n})$ directed by $\alpha$.
This fact motivates referring to $\alpha$ as the \textit{directing random measure} of the c.i.d. sequence.
Further asymptotic results and uniform limit theorems for c.i.d. sequences were given, among others, by \cite{BertiPratelliRigo2004LimThForIID, 
BertiPratelliRigo2012LimThforEmpirProcBasedOnDepData}.

\section{Partially c.i.d. sequences} \label{sec:pcid}

Our aim is to extend the notion of c.i.d. sequences to arrays of random variables $[X_{n,i}]_{n \geq 1; i \in I}$, 
 where $I$ is a finite or countable set. We introduce a notion of {\em partially c.i.d.} sequences, that can be regarded as a weaker form of partial exchangeability, in the same sense that the c.i.d. property is a weaker form of exchangeability.

An array of random variables $[X_{n,i}]$ 
is partially exchangeable (in the sense of de Finetti) if its probability law is invariant under separate finite permutations of the columns,  that is
\begin{equation*} 
[X_{n,i}]\eqdist [X_{ \pi_i(n), i}]
\end{equation*}
for all finite permutations $\pi_i$, $i\in I$.  Marginally, the columns $(X_{n,i})_{n\geq 1}$ are exchangeable, but a family of exchangeable sequences is not necessarily partially exchangeable. For example, two exchangeable sequences $(X_n)$ and  $(Y_n)$, where $X_n=Y_n$ for all $n$, are not partially exchangeable, unless their distributions are degenerate. Similarly, an interesting notion of partially c.i.d. dependence should  require appropriate conditions for a family of c.i.d. sequences to be partially c.i.d.
We proceed along the lines of Section \ref{sec:cid} and, as a first step, we extend the notion of spreadability to families of sequences. Let us introduce the notion of \textit{partial spreadability} as invariance of the joint probability law under separate selection of subsequences along the columns.
\begin{defn}
An array $[X_{n,i}]_{n \geq 1; i\in I}$ of random variables is {\em partially spreadable} if 
\begin{equation} \label{eq:pspread}
[X_{n, i}] \eqdist [X_{k_n^{(i)}\!\!\!,\,i}],   \quad \mbox{for every }   k^{(i)}_1 < k^{(i)}_2 < \cdots, \;  i \in I.
\end{equation}
\end{defn}
Clearly, a partially exchangeable array satisfies (\ref{eq:pspread}). The reverse implication is also true.
\begin{prop}
A partially spreadable array $[X_{n,i}]$ is partially exchangeable.
\end{prop}
\noindent {\sc Proof}.
Let $J$ be a finite subset of $I$ and $j$ be a fixed element in $J$. Let $Z=[X_{n,i}]_{n\geq 1; i\in J\setminus\{j\}}$,
and, for every $n\geq 1$, let $U_n=(Z,X_{n,j})$. By (\ref{eq:pspread}), the sequence $(U_n)$ is spreadable. By the results in \cite{ryllNad1957},  $(U_n)$ is exchangeable. In other words, $(X_{n,j})_{n\geq 1}$ is exchangeable over $Z$ (\cite{Aldous85}, Section 3). Since this holds for every finite $J\subset I$ and for every $j\in J$, the array $[X_{n,i}]$ is partially exchangeable (\cite{Aldous85}, Proposition 3.8). \\
\noindent $\square$

Let us now regard the rows of $[X_{n,i}]$ as the values of a process $(\bX_n^I)_{n \geq 1}$,  where $\bX_n^I=(X_{n,i}, i \in I)$.
If the process $(\bX_n^I)$ is stationary, then a weaker form of spreadability is sufficient for partial exchangeability.
The following result extends Proposition \ref{prop:kallenberg}.
\begin{thm} \label{prop:p-kallenberg}   
A  stationary process  $(\bX_n^I)$ that satisfies
\begin{equation}    \label{eq:pcid}
( {\bf X}^{J}_{1:n}, \bX^{J\setminus\{j\}}_{n+1}, X_{n+1,j})
\eqdist ( {\bf X}^{J}_{1:n}, \bX^{J\setminus\{j\}}_{n+1}, X_{n+2,j})\quad \mbox{for every } n\geq 0,j\in J\mbox{ and every finite }J\subset I,
\end{equation}
is partially exchangeable, where (\ref{eq:pcid}) reads $(\bX^{J\setminus\{j\}}_{1}, X_{1,j})
 \eqdist ( \bX^{J\setminus\{j\}}_{1}, X_{2,j}) $ when $n=0$.
\end{thm}


\noindent {\sc Proof}.
Since $(\bX_n^I)$ is stationary, we can embed it in a doubly infinite sequence $(\bX_{n}^I)_{-\infty<n<\infty}$.
Marginally, each sequence $(X_{n,i})$ is stationary and c.i.d., therefore, by Proposition \ref{prop:kallenberg}, exchangeable. Let us denote by $\alpha_i$ its directing random measure.
Let $J$ be a finite subset of $I$, $j$ a fixed element of $J$ and $m$ a fixed positive integer.
 For every $n\geq 1$, let  $Z_n=(\bX_{1:m}^J,\bX_{m+1}^{J\setminus\{j\}},X_{m+n,j})$. Since $(Z_n)$ is stationary and c.i.d., it is exchangeable. 
 In other words, $(X_{m+n,j})_{n\geq 1}$ is exchangeable over $(\bX_{1:m}^J,\bX_{m+1}^{J\setminus\{j\}})$
 (\cite{Aldous85} Section 3).
 By Proposition 3.8 in \cite{Aldous85}, $(X_{m+1,j},X_{m+2,j},\dots)$ are conditionally i.i.d. with common law $\alpha_j$,
 given $(\alpha_j, (\bX_{1:m}^J,\bX_{m+1}^{J\setminus\{j\}}))$.
 Since $(\bX_{n}^I)$ is stationary, for every $k\in\mathbb Z$, 
 $X_{k,j}$ is conditionally independent of $(X_{h,i}:(h,i)\neq (k,j),k-n\leq h\leq k,i\in J)$, 
 given $\alpha_j$. Therefore,  $X_{k,j}$ is conditionally independent of 
 $(X_{h,i}:(h,i)\neq (k,j),h\leq k,i\in J)$, given $\alpha_j$. 
 Since $\alpha_i=\lim_{n\rightarrow\infty}\sum_{s=1}^n\delta_{X_{k-s,i}}/n$ ($\mathbb  P\mbox{-a.s.}$), $X_{k,j}$ is also conditionally independent of $(\alpha_i:i\neq j)$, given $\alpha_j$. 
 Thus, for every $n$
 $$
 \mathbb P[\cap_{i\in J}\cap_{k=1}^n (X_{k,i}\in A_{k,i})\mid \alpha_i,i\in J]=\prod_{i\in J}\prod_{k=1}^n \alpha_i( A_{k,i}), \quad \mathbb P\mbox{-a.s.}
 $$
 and $(\bX_n^I)$ is partially exchangeable.\\
 \noindent $\square$

Condition (\ref{eq:pcid}) implies that, for any $j \in I$, future values of $X_{n,j}$ are conditionally identically distributed, given the past observations $(\bX_{1:n}^I)$ and the concomitant values $(\bX_{n+1}^{I\setminus\{j\}})$ of the other variables.
More generally, given a filtration $\mathcal{G}=\left(\mathcal{G}_{n}\right)_{n\geqslant 0}$, we introduce the following definition.
\begin{defn}
A sequence $(\bX^I_n)$ is said to be partially conditionally identically distributed with respect to a filtration $\mathcal{G}$ (briefly, partially $\mathcal{G}$-c.i.d.), if it is adapted to $\mathcal G$ and, for every $j \in I$, 
\begin{equation}\label{eq:PCID1}
\begin{aligned}
&\mathbb E[f(X_{n+1, j})  \mid \calG_n^j]=\mathbb E[f(X_{n+k, j}) \mid \calG_n^j]\\
&\mbox{ for all } k \geq 1, n \geq 0 \mbox{ and all bounded measurable functions } f:\mathbb X\rightarrow \mathbb R,
\end{aligned}
\end{equation}
where $\mathcal G_n^j=\mathcal G_n \vee \sigma(X_{n+1,i} : i \neq j)$.
\end{defn}
When $\mathcal G$ is the natural filtration of $(\bX_n^I)$, with $\mathcal G_0=\{\emptyset,\Omega\}$, the process is said to be naturally partially c.i.d. or, simply, partially c.i.d.
Clearly, in this case condition  \eqref{eq:PCID1} reduces to \eqref{eq:pcid}.
It is easy to show that, if $(\bX_n^I)$ is partially $\mathcal{G}$-c.i.d., then it is also partially c.i.d. with respect to any coarser filtration to which it is adapted; in particular,  it is partially c.i.d.
A partially exchangeable array satisfies \eqref{eq:pcid}, thus it is partially c.i.d.

Condition \eqref{eq:PCID1}  means that each sequence $(X_{n,j})_{n \geq 1}$ is c.i.d. with respect to the filtration $\mathcal G^j=(\mathcal G_n^j)_{n\geq 0}$. This allows us to formulate equivalent conditions (extending Proposition \ref{prop:kallenberg2}).
\begin{prop}
The following properties (i), (ii) and (iii) are equivalent.
\begin{itemize}
\item[(i)] The sequence $(\bX_n^I)$ is partially $\calG$-c.i.d.;
\item[(ii)] For any $j \in I$ and every 
bounded measurable function $f:\mathbb X\rightarrow \mathbb R$,
$ \left(\mathbb E\left[f(X_{n+1,j}) \mid \calG_n^j\right]\right)_{n\geqslant0}$ is a $\calG^j$-martingale;
\item[(iii)]
For every $n\geq 0$, and $j\in I$,
	$(\bX_{n+1}^{I\setminus \{j\}},X_{\tau+1, j})
	\eqdist
	(\bX_{n+1}^{I\setminus \{j\}},X_{n+1, j}),$
for all finite $\mathcal G$-stopping times $\tau$ satisfying $\tau \geq n$.
\end{itemize}
\end{prop}

\noindent \textsc{Proof.} 
Conditions (i) and (ii) are equivalent, by Proposition \ref{prop:kallenberg2}. To prove that (ii) implies (iii), fix $n\geq 0$ and $j\in I$ and let $\tau \geq n$ be a finite 
$\mathcal G$-stopping time. Let $f:\mathbb X\rightarrow \mathbb R$ be a bounded measurable function and let $B\in \mathcal X^{I\setminus\{j\}}$.  Since $(\mathbb E[f(X_{n+k+1,j})\mathbbm 1_{\{\bX_{n+1}^{I\setminus \{j\}}\in B\}}\mid \mathcal G_{n+k}])_{k\geq 0}$ is a bounded martingale and $\tau-n$ is a finite stopping time,  with respect to the filtration $(\mathcal G_{n+k})_{k\geq 0}$, by Doob optional stopping theorem,
 $\mathbb E[f(X_{\tau+1,j})\mathbbm 1_{\{\bX_{n+1}^{I\setminus \{j\}}\in B\}}]=\mathbb E[f(X_{n+1,j})\mathbbm 1_{\{\bX_{n+1}^{I\setminus \{j\}}\in B\}}]$, which implies (iii).
To prove that   (iii) implies (i),
let us fix $B\in\mathcal X$, $C\in\mathcal X^{I\setminus \{j\}}$, $A\in\mathcal{G}_{n}$
and let $\tau=(n+1) \mathbbm{1}_A+n\mathbbm{1}_{A^c}$. Then, $\tau$ is a finite $\mathcal{G}$-stopping time and
$$
\begin{aligned}
\mathbb{P}\!\left[(X_{n+2,j}\in B)\!\cap\! (\bX^{I\setminus\{j\}}_{n+1}\!\in \! C)\!\cap\! A\right]
\!&=\!\mathbb{P}\!\left[(X_{\tau+1,j}\in B)\!\cap\! (\bX^{I\setminus\{j\}}_{n+1}\!\in\! C)\right]\!-\!\mathbb{P}\!\left[(X_{\tau+1,j}\!\in\! B)\!\cap \!(\!\bX^{I\setminus\{j\}}_{n+1}\!\in\! C\!)\!\cap\!(\tau=n)\right]\\
&=\mathbb{P}\left[(X_{n+1,j}\in B)\cap (\bX^{I\setminus\{j\}}_{n+1}\in C)\cap A\right]. 
\end{aligned}
$$
\noindent $\square$

Theorem \ref{prop:p-kallenberg} shows that stationarity is a sufficient condition for a partially c.i.d. process to be partially exchangeable. 
The following proposition gives a sufficient condition in terms of the predictive distributions.

\begin{prop} \label{prop:predCharacterization}
	Let $(\bX_n^I)$ be a partially c.i.d. process. Suppose that, for every finite $J\subset I$ and every $j\in J$,  $n\geq 0$, permutations  $(\pi_i)_{i\in J}$, we have
		\begin{equation}\label{eq:repres_p_exch}\begin{aligned}
			\mathbb{P}\left[X_{n+1,j}\in\cdot\mid \cap_{(m,i)\in H_{n,j}}(X_{m,i}\in A_{m,i})\right]		=\mathbb{P}\left[X_{n+1,j}\in\cdot\mid \cap_{(m,i)\in H_{n,j}}(X_{\pi_i(m),i}\in A_{m,i})\right] , \\
			\mbox{ where } H_{n,j}=\{(m,i)\neq (n+1,j):m\leq n+1, i\in J\},
		\end{aligned}
	\end{equation}
		 for all $A_{m,i}$  such that the conditioning sets have positive probability. Then $(\bX_n^I)$ is partially exchangeable.
		\end{prop}
\noindent \textsc{Proof.}
By Theorem \ref{prop:p-kallenberg}, we only need to prove that $(\bX_{n}^I)$ is a stationary process, which is true if $(\bX_{n}^J)$ is stationary for every finite $J\subset I$. 
Let $J$ be a finite subset of $I$.
Without loss of generality, we can assume that $J=\{1,2,\dots,d\}$, for some $d\in \mathbb N$. Let us define the following ordering in $\mathbb N\times J:$ 
\begin{equation}
	\label{eq:order}
	(1,1)\prec(1,2)\prec\dots \prec(1,d)\prec(2,1)\prec(2,2)<\dots \prec(2,d)\prec(3,1)\prec\dots\prec(3,d)\prec\dots.
\end{equation}
For every $N\in\mathbb N$, let us denote by $K_N\subset \mathbb N\times J$ the set of the first $N$ pairs in (\ref{eq:order}). We prove that
\begin{equation}
	\label{eq:stat}
	\mathbb P[\cap_{(m,i)\in K_N}(X_{m,i}\in A_{m,i})]=\mathbb P[\cap_{(m,i)\in K_N}(X_{m+1,i}\in A_{m,i})]\quad \mbox{ for every } [A_{m,i}]_{(m,i)\in K_N}\in\mathcal X^N
\end{equation}
for every $N$, which implies stationarity of $(\bX_{n}^J)$. The proof is done by induction on $N$. The property is trivially true for $N=1$, since $(X_{n,1})_{n\geq 1}$ is c.i.d. 
Suppose that (\ref{eq:stat}) holds for some $N\in \mathbb N$ and let $A_{m,i}$, with $(m,i)\in K_{N+1}$, be fixed subsets in $\mathcal X$. Let $(n,j)$ denote the largest element of $K_N$ and let $(n',j')$ denote its subsequent element in (\ref{eq:order}). Then $n'\geq n$ and $j'\neq j$.
If $\mathbb P[\cap_{(m,i)\in K_N}(X_{m,i}\in A_{m,i}) ]=0$, then, by the induction hypothesis, $\mathbb P[\cap_{(m,i)\in K_N}(X_{m+1,i}\in A_{m,i}) ]=0$ and (\ref{eq:stat}) holds with $N+1$ in the place of $N$, since both terms of the equality are equal to zero. 
Conversely, suppose that $\mathbb P[\cap_{(m,i)\in K_N}(X_{m,i}\in A_{m,i}) ]>0$. By the induction hypothesis $\mathbb P[\cap_{(m,i)\in K_N}(X_{m+1,i}\in A_{m,i}) ]>0$.
 By the p-c.i.d. condition and the induction hypothesis,
\begin{align*}
	\mathbb P\left[\cap_{(m,i)\in K_{N+1}}(X_{m,i}\in A_{m,i})\right]&=\mathbb P\left[X_{n',j'}\in A_{n',j'}\mid \cap_{(m,i)\in K_N}(X_{m,i}\in A_{m,i})\right]
	\mathbb P\!\left[\cap_{(m,i)\in K_N}(X_{m,i}\in A_{m,i}) \right]\\
	&=\mathbb P[X_{n'+1,j'}\in A_{n',j'}\mid \!\cap_{(m,i)\in K_N}(X_{m,i}\in\! A_{m,i})]
	\,\mathbb P\!\left[\cap_{(m,i)\in\! K_N}(X_{m+1,i}\in \!A_{m,i}) \right]\!.
\end{align*}
On the other hand, since $(n'+1,j')$ is the largest element in $K_{N+d+1}$, by (\ref{eq:repres_p_exch}),
\begin{align*}
	&\mathbb P[X_{n'+1,j'}\in A_{n',j'}\mid \cap_{(m,i)\in K_N}(X_{m,i}\in A_{m,i})]
	=\mathbb P[X_{n'+1,j'}\in A_{n',j'}\mid \cap_{(m,i)\in K_N}(X_{m+1,i}\in A_{m,i})].
\end{align*}
Hence (\ref{eq:stat}) holds for $N+1$ in the place of $N$, and the thesis follows. 

\noindent $\square$

\begin{cor} \label{prop:repres_thm}
A c.i.d. sequence $(X_n)$ is exchangeable if and only if, for every $n\geq 1$ and permutation $\pi$,
	\begin{equation*} 
	\mathbb{P}\left[X_{n+1}\in \cdot\mid X_1\in A_1,\dots,X_n\in A_n\right]=\mathbb{P}\left[X_{n+1}\in \cdot\mid X_{\pi(1)}\in A_1,\dots,X_{\pi(n)}\in A_n\right],
		\end{equation*}
		for all $A_1,\dots, A_n$ such that the conditioning sets have positive probability.
\end{cor} 
Intuitively, a c.i.d. sequence is exchangeable if the order of the past observations is \textit{totally irrelevant} for predictions based on whatever information about the past. A different set of predictive conditions that characterize an  exchangeable probability law is given by \cite{fortiniLadelliRegazzini2000}.

 \begin{exa}  \label{ex:interactingRRprocesses}  \rm {\bf [Interacting randomly reinforced processes]}.
 As discussed in the Introduction, processes with reinforcement are a powerful way to generate exchangeable sequences. However, even natural extensions to systems of dependent reinforced processes break stationarity and are not partially exchangeable. In fact,
they are partially c.i.d. under quite natural assumptions on the updating rule.
Consider a family of sequences $(X_{n,i})_{n \geq 1}$, $ i \in I$, that evolve according to the following randomly reinforced scheme. Let $w_{0,i}$ be known positive scalars and $W_{1,i}, W_{2,i}, \ldots$ be positive random variables, for every $i\in I$. Let  $\mathcal G_n=\sigma(\bX_{1:n}^I,\mathbf W_{1:n}^I)$.
Assume that $X_{1,i} \sim \nu_i(\cdot)$ for a given distribution $\nu_i$, and for any $n \geq 1$
\begin{equation} \label{eq:pred-RRprocess}
Q_{n,i}(\cdot) := {\mathbb P}[X_{n+1,i} \in \cdot \mid \calG_n] =
\frac{w_{0,i} \nu_i(\cdot) + \sum_{k=1}^n W_{k,i} \delta_{X_{k,i}}(\cdot)}
{w_{0,i} + \sum_{k=1}^n W_{k,i}}, \quad i \in I.
\end{equation}
The generative rule (\ref{eq:pred-RRprocess}) produces a system of reinforced processes. For binary random variables, the interacting randomly reinforced urns by \cite{paganoniSecchi2004} are examples of the construction  (\ref{eq:pred-RRprocess}).
Other important  reinforced schemes are special cases.
If $W_{n,i}=1$ for all $n$ and $i$, (\ref{eq:pred-RRprocess}) describes independent \Polya sequences \citep{blackwellMacqueen1973}, that is, independent, internally exchangeable sequences $(X_{n,i})_{n \geq 1}$ driven by a Dirichlet process with base measure $w_{0,i} \, \nu_i$.  
For time-varying or random weights $W_{n,i}$, marginal exchangeability is generally lost; however, such extensions are of interest in many applications. For a single sequence (i.e. $I=\{1\}$),
 a randomly reinforced scheme of the kind (\ref{eq:pred-RRprocess}) is  
considered by \cite{Bassetti2010CidSSS} as a (c.i.d.) generalization of  exchangeable species sampling models; further developments are given in \cite{airoldiCosta2014}. The predictive system (\ref{eq:pred-RRprocess}) envisages extensions to multiple populations. 

Indeed, the interest in the generative rule (\ref{eq:pred-RRprocess}) is that one can introduce interaction across the sequences, through the weights $W_{n,i}$, which are generally stochastically dependent;
in particular, each $W_{n,i}$ may be a function of the observed values of the other sequences. 
Nonetheless, it is natural to assume that the populations are updated independently, that is,

(a) conditionally on $\calG_n$, the $X_{n+1,i}$, $i \in I$, are independent.\\
If, moreover,

(b) the random weight $W_{n,i}$ is conditionally independent of $X_{n,i}$, given $\mathcal G_{n-1}$, for any $n \geq 1, i \in I$, \\
then $(\bX_n^I)$
is partially $\calG$-c.i.d.
To see this, write the predictive distribution (\ref{eq:pred-RRprocess}) recursively:
	\begin{equation} \label{eq:recursive}
	Q_{n,i}(\cdot)= A_{n,i} \delta_{X_{n,i}}(\cdot) + (1-A_{n,i})  Q_{n-1,i}(\cdot), 
	\end{equation}
	where $A_{n,i}= W_{n,i}/(w_0+\sum_{k=1}^n W_{k,i})$. 
By (\ref{eq:recursive}) and assumptions (a) and (b), for every $n\geq 0$,
		\begin{eqnarray*}
		\mathbb E[Q_{n+1,i} (\cdot) \mid \calG_n^i\vee \sigma(W_{n+1,i})]
				=\mathbb  E \left[A_{n+1,i}\delta_{X_{n+1,i}}(\cdot)+(1-A_{n+1,i})Q_{n,i} (\cdot) \mid \calG_n^i\vee \sigma(W_{n+1,i})\right]=Q_{n,i}(\cdot).
	\end{eqnarray*}
	Thus, $[X_{n,i}]$ is partially $\calG$-c.i.d. Notice that (\ref{eq:pred-RRprocess}) and (\ref{eq:recursive}) are equivalent. Indeed, (\ref{eq:pred-RRprocess}) can be retrieved from (\ref{eq:recursive}) by defining, recursively, $w_{0,i}=1$ and $W_{n,i}=(w_0+\sum_{k=1}^{n-1}W_{k,i})A_{n,i}/(1-A_{n,i})$.

\medskip
Forms of interaction in the partially c.i.d. system (\ref{eq:pred-RRprocess}) may include common weights: $W_{n,i}= W_n$ for any $i \in I$ and $n \geq 1$; interaction through the weights: $W_{n,i}=w_{n,i}(X_{k,j}:k\leq n,j\neq i)$; interaction through common (observable or latent) variables:
$W_{n,i}= w_{n,i}(Z_{k}:k\leq n)$ for any $i$. In the latter, $\calG_n= \sigma( \bX_{1:n}^I, Z_{1:n})$, and one further assumes that $Z_{n+1}$ is conditionally independent of $\bX_{n+1}^I$, given $\calG_n$. A further extension is to let  $W_{n,i}= w_{n,i}(Z_{1:n}, \bX_{1:n}^{I\setminus\{i\}})$.
\end{exa}

The following is an example of the construction (\ref{eq:pred-RRprocess}), where no latent  variables are introduced, so that the predictive distribution given the natural filtration has a simple form.

\begin{exa}
	\label{ex:uniform}\rm
Consider the system (\ref{eq:pred-RRprocess}) where $I=\{1,2\}$, the marginal law $\nu_i$ is the uniform distribution on the unit interval and the random weights $W_{n,i}$ are defined recursively through  (\ref{eq:recursive}), by letting $A_{n,i}= \beta_n X_{n,j}$ for $j \neq i$, where $(\beta_n)$ is a fixed sequence of positive numbers, with $\beta_1=1$ and $\beta_n\leq 1$. 
These assumptions could be relaxed, allowing for more than two sequences and for different marginal laws.
The independence assumptions (a) and (b) in Example \ref{ex:interactingRRprocesses} are satisfied, thus 
$(\bX_n^I)$
is partially-c.i.d.
The sequence of parameters $(\beta_n)$ allows to tune the relative weights of past observations in the predictive distributions. For example, if $\beta_n=1$ for every $n$, then the expected weight of $X_{n,i}$ in the prediction of $X_{n+1,i}$ is $1/2$, and observations that are far away in the past have negligible weight; conversely, if $\beta_n=2/(n+1)$, then the expected  weight of $X_{n,i}$ is $1/(n+1)$ and all past observations have roughly the same weight.

The interaction between the two sequences induces a positive correlation in the joint law of $(X_{n,1},X_{n,2})$, for every $n$. Indeed,
$\mathbb E[(X_{2,1}-1/2)(X_{2,2}-1/2)]
=\mathbb E[(X_{1,1}-1/2)^2(X_{1,2}-1/2)^2]=1/144,$
so that $Corr(X_{2,1},X_{2,2})=1/12$. Furthermore, denoting by $\mathcal G$ the natural filtration of $(\bX_n^I)$,
\[
\begin{aligned}
&\mathbb E[(X_{n+1,1}-1/2)(X_{n+1,2}-1/2)]-\mathbb E[(X_{n,1}-1/2)(X_{n,2}-1/2)]
=\beta_n^2\mathbb E[\mathbb V[X_{n,1}\mid\mathcal G_{n-1}]\mathbb V[X_{n,2}\mid\mathcal G_{n-1}]]\geq 0.
\end{aligned}
\]
Since the standard deviations are constant, the correlation between $X_{n,1}$ and $X_{n,2}$ is an increasing function of $n$ and is, therefore, strictly positive for every $n>1$. 
\end{exa}

\section{Asymptotic partial exchangeability} \label{sec:partialExch}
As seen in Section \ref{sec:cid}, c.i.d. sequences are asymptotically exchangeable. Here, we prove that partially c.i.d. processes are asymptotically partially exchangeable.

Clearly, asymptotic partial exchangeability does not follow from marginal properties. If an array  $[X_{n,i}]$ is partially $\calG$-c.i.d., then
each of its columns is $\calG$-c.i.d.
and the marginal predictive distributions $\mathbb{P}[X_{n+1, i}\in \cdot \mid\mathcal{G}_n]$  converge to the {\em marginal directing measure} $\alpha_i$ on $\mathbb{X}$. 
This implies that the vector $(\mathbb{P}[X_{n+1, i}\in \cdot \mid\mathcal{G}_n], i \in I)$ converges to the vector of random measures $(\alpha_i, i \in I)$.
However, these properties are not  sufficient to ensure the asymptotic conditional independence required for partial exchangeability.
Let us recall that, by the representation theorem for partially exchangeable sequences (see  \cite{Aldous85}, Section 3), an array $[X_{n,i}]$ is partially exchangeable if and only if there exist random probability measures $(\alpha_i, i \in I)$ such that, conditionally on $(\alpha_i, i \in I)$, the  sequences $(X_{n,i})$ are independent, with $X_{n,i} \mid (\alpha_j, j \in I) \indsim \alpha_i$, $i \in I$. 

In order to prove the asymptotic partial exchangeability of partially c.i.d. arrays, we
show that the {\em joint} predictive distribution  $\mathbb P[\bX_{n+1}^I \in \cdot \mid \calG_n]$
converges a.s. (which is not an immediate consequence of the marginal c.i.d. property, because the sequence $(\bX_n^I)$ is not c.i.d.); and that it converges to the product random probability measure $\mathlarger{\mathlarger{\mathlarger{\times}}}_{i \in I} \alpha_i$.
This implies asymptotic partial exchangeability (Theorem \ref{th:pexc}).

Moreover, we show that, as for partially exchangeable sequences, the empirical distributions converge a.s. to the same limit as the predictive distributions.  Here, convergence is meant as weak convergence. However, with a little more effort, we can prove that the predictive and the empirical distributions converge \textit{point-wise}. This is the content of  Proposition \ref{prop:strong}. The point-wise convergence will also become useful in Section \ref{sec:SLLN}. The basic step for the proof is given by the following Lemma.

\begin{lemma}
	\label{lemma2}
	Let $(\bX_n^I)$
		be a partially $\mathcal{G}$-c.i.d. process. Then, for every finite $J \subset I$ and every  bounded and measurable functions $f_i$, $i\in J$,
	\begin{equation}
		\label{eq:pcidconv2}
		\mathbb{E}\left[\prod_{i\in J} f_i(X_{n+1,i})\mid\mathcal{G}_n\right]\rightarrow \prod_{i\in J} \int f_id\alpha_i\quad \mathbb P\mbox{-a.s. and in }L^1,
	\end{equation}
	and
	\begin{equation}
		\label{eq:conv_emp}
		\frac{1}{n}\sum_{k=1}^n\prod_{i\in J} f_i(X_{k,i})
				\rightarrow \; \prod_{i\in J}\int f_id\alpha_i
		\quad \mathbb P\mbox{-a.s. and in }L^1,
	\end{equation}
	as $n\rightarrow \infty$, where $\alpha_i$ is the marginal directing measure of $(X_{n,i})_{n \geq 1}$.	
\end{lemma}

\noindent\textsc{Proof.}
Since the $f_i$ are bounded, $L^1$ convergence is a consequence of the a.s. convergence.\\
The proof of (\ref{eq:pcidconv2}) proceeds by induction on $k=|J|$, the cardinality of $J$. For $k=1$, let $J=\{i\}$. Since $(X_{n,i})$ is $\mathcal G$-c.i.d., (\ref{eq:pcidconv2}) holds for indicators  
\citep{BertiPratelliRigo2004LimThForIID} and extends by linearity to simple functions. 
Now let $f$ be bounded and measurable. Then, $\mathbb E[f(X_{n+1,i})\mid \mathcal G_n]$ is a uniformly integrable martingale. Let $V_f$ be its $\mathbb P$-a.s. limit. To prove that $V_f=\int fd\alpha_i$, it is sufficient to approximate $f$ from below with a sequence $(g_m)_{m\geq 1}$ of simple functions satisfying $f-g_m<1/m$. Then, $\mathbb P$-a.s.  $\int fd\alpha_i-\int g_m d\alpha_i<1/m$ and $\mathbb E[f(X_{n+1,i})\mid \mathcal G_n] -\mathbb E[g_m(X_{n+1,i})\mid \mathcal G_n] <1/m$, which implies that $V_f-V_{g_m}<1/m$. Since, for every $m$, $V_{g_m}=\int g_m d\alpha_i$, then $V_f=\int fd\alpha_i$, $\mathbb P$-a.s. Thus, (\ref{eq:pcidconv2}) holds for $|J|=1$. Let us now prove that, if (\ref{eq:pcidconv2}) holds for some $J$, then it also holds for $J\cup\{j\}$, where $j$ is any element of $I\setminus J$.
 By the definition of partially c.i.d. array,
\[
\mathbb{E}[\prod_{i\in J\cup\{j\}} f_i(X_{n+1,i})\mid\mathcal{G}_n]
=\mathbb{E}[\prod_{i\in J} f_i(X_{n+1,i}) \int f_{j}d\alpha_{j}\mid\mathcal{G}_n]
=A_n+\mathbb{E}[\prod_{i\in J} f_i(X_{n+1,i})\mid\mathcal{G}_n]\,\mathbb E[f_{j}(X_{n+1,j})\mid\mathcal{G}_n]\\
\]
where
$$
|A_n| \leq \mathbb E\left[\left.\prod_{i\in J}\left| f_i(X_{n+1,i})\right|   \left| \int f_{j}d\alpha_{j}-\mathbb E[\int f_{j}d\alpha_{j}\mid \mathcal{G}_n]\right| \right| \mathcal{G}_n\right]\rightarrow 0 \quad\mathbb P\mbox{-a.s.}\quad \mbox{as }n\rightarrow\infty,
$$
since $ \prod_{i\in J} |f_i(X_{n+1,i})|$ is bounded and $|\int f_{j}d\alpha_{j}-\mathbb E[\int f_{j}d\alpha_{j}\mid \mathcal{G}_n]|=|\int f_{j}d\alpha_{j}-\mathbb E[f_{j}(X_{n+1,j})\mid \mathcal{G}_n]|$ is bounded and converges to zero a.s. 
Using the induction hypothesis, we obtain (\ref{eq:pcidconv2}) for $J\cup\{j\}$ in the place of $J$.
\\
Let us now prove (\ref{eq:conv_emp}). Since the left hand side in (\ref{eq:conv_emp}) is uniformly bounded, the $L^1$ convergence follows from almost sure convergence. Let $Y_k=\prod_{i\in J} f_i(X_{k,i})$ and
\[
M_n=\sum_{k=0}^{n-1}\frac{Y_{k+1}-\mathbb E[Y_{k+1}\mid \mathcal G_k]}{k+1}.
\]
Then $(M_n)_{n\geq 1}$ is a martingale and satisfies
$\sup_n \mathbb E[M_n^2]<\infty$. 
Hence, $M_n$ converges almost surely. By the Kronecker Lemma,
\[
\frac{1}{n}\sum_{k=0}^{n-1}\left(Y_{k+1}-\mathbb E[Y_{k+1}\mid\mathcal G_k]\right)\rightarrow 0 , \quad \mathbb P\mbox{-a.s.}
\]
By (\ref{eq:pcidconv2}),   $
\sum_{k=0}^{n-1}\mathbb E[Y_{k+1}\mid\mathcal G_k]/n\rightarrow \prod_{i\in I}fd\alpha_i$, $\mathbb P$-a.s.
Hence, (\ref{eq:conv_emp}) holds. \\
\noindent $\square$

\begin{prop}
	\label{prop:strong}
	Let $(\bX_n^I)$ be a partially $\calG$-c.i.d. process,
	with marginal directing measures $(\alpha_i: i \in I)$. Then, there exists $N\in \mathcal F$ with $\mathbb P[N]=0$ such that, for every $\omega\in N^c$,
	\begin{equation}
		\label{eq:strong}
		\lim_{n\rightarrow\infty}\mathbb P[ \bX_{n+1}^I \in B \mid \mathcal G_n](\omega)=\lim_{n\rightarrow\infty}\frac{1}{n}\sum_{k=1}^n\delta_{\bX_k^I}(B)(\omega) =\alpha(B)(\omega)   \quad \mbox{ for every }B \in \mathcal X^I,
	\end{equation}
		and
	\begin{equation}
		\label{eq:strong2}
		\lim_{n\rightarrow\infty}\mathbb E[f(\bX_{n+1}^I)\mid\mathcal G_n](\omega)=\lim_{n\rightarrow\infty}\frac{1}{n}\sum_{k=1}^n f(\bX_k^I)(\omega)= \int fd\alpha(\omega )\quad \mbox{for every bounded measurable} f ,
	\end{equation}
	where $\alpha=\mathlarger{\mathlarger{\mathlarger{\times}}}_{i\in I}\alpha_i$.
	\end{prop}
\noindent\textsc{Proof.}
 Let $P_{n+1}(\cdot):=\mathbb P[\bX_{n+1}^I\in \cdot\mid\mathcal G_n]$, $\hat P_n(\cdot):=n^{-1}\sum_{k=1}^n\delta_{\bX_k^I}(\cdot)$,
 $E_{n+1}(f):=\mathbb E[f(\bX_{n+1}^I)\mid \mathcal G_n]$ and $\hat E_n(f):=n^{-1}\sum_{k=1}^n f(\bX_k^I)$.\\
Since $\mathcal X$ is countably generated, by Lemma \ref{lemma2}, there exists $N\in\mathcal F$, with $\mathbb P[N]=0$, and a class $\mathcal C\subset \mathcal X$ generating $\mathcal X$, closed under complements, finite intersections and disjoint finite unions, such that, 
for all $B\in \mathcal S=\{\times_{i\in J}A_i\times \mathbb X^{I\setminus J}: J\mbox{ finite, }A_i\in\mathcal C, i\in J \}$, we have
\begin{equation}
\label{eq:lambda}
\lim_{n\rightarrow\infty} P_{n+1}(B)(\omega)=\lim_{n\rightarrow\infty}\hat P_n(B)(\omega)= \alpha(B)(\omega)\quad \mbox{for every } \omega\in N^c.
\end{equation}
%
%
%
Property (\ref{eq:lambda}) extends easily to countable disjoint unions of events in $\mathcal S$ and, using Dynkin Lemma, to every  $B\in \mathcal X^I$.
	This proves (\ref{eq:strong}).
\\	We now prove (\ref{eq:strong2}). If $f$ is a simple function, by (\ref{eq:strong}), $\lim_{n\rightarrow\infty}E_{n+1}(f)=\lim_{n\rightarrow\infty} \hat E_n(f)=\int fd\alpha$ for every $\omega\in N^c$. If $f$ is non negative,  for every $\epsilon>0$, there exists a simple function $f^{(\epsilon)}$ such that  $0\leq f-f^{(\epsilon)} <\epsilon$. Hence,  for every $\omega\in N^c$,
\[
\liminf_{n\rightarrow\infty}  E_{n+1}(f)(\omega)
\geq \liminf_{n\rightarrow\infty} E_{n+1}(f^{(\epsilon)})(\omega)
\geq \int f^{(\epsilon)}d\alpha(\omega)
\geq \int fd\alpha(\omega)-\epsilon.
\]
Since the above inequalities hold for every $\epsilon >0$, $\liminf E_{n+1}(f)(\omega)\geq \int fd\alpha(\omega)$.
 If $M=\sup f$, then $M-f$ is bounded and non negative; hence
	$\limsup E_{n+1}(f)(\omega)\leq M-\int(M-f)d\alpha(\omega)\leq \int fd\alpha(\omega)$.
			Thus the limit $\lim_n\mathbb E_{n+1}(f)(\omega)
			=\int f d\alpha(\omega)$ exists. By linearity, the property extends to functions $f$ taking positive and negative values.
By repeating the same reasoning, with $\hat E_n$ replacing $E_{n+1}$, we obtain the thesis. \\
\noindent $\square$

From Proposition \ref{prop:strong}, it follows that the predictive distributions
$\mathbb P[\bX_{n+1}^I \in \cdot \mid\mathcal G_n]$ converge weakly to the product random probability measure $\alpha=\mathlarger{\mathlarger{\mathlarger{\times}}}_{i\in I} \alpha_i$, $\mathbb P$-a.s.
This implies that $(\bX_n^I)$ is asymptotically exchangeable  with directing measure $\alpha$ \citep[Lemma 8.2]{Aldous85}.
Being $\alpha$ a product measure,  $(\bX_n^I)$ is asymptotically partially exchangeable. This proves the following 
\begin{thm}
	\label{th:pexc}
Let  $(\bX_n^I)$ be a partially $\mathcal{G}$-c.i.d. process with marginal directing measures $(\alpha_i,i\in I)$. 
Then, $(\bX_n^I)$ is asymptotically partially exchangeable, and the partially exchangeable limit law has directing measures $(\alpha_i,i\in I)$.
\end{thm}

\medskip
Theorem \ref{th:pexc} justifies referring to  $(\alpha_i, i \in I)$ as  the directing random measures of the partially c.i.d. sequence $(\bX_n^I)$.
From Theorem \ref{th:pexc} 
it follows that a stationary partially c.i.d. sequence is partially exchangeable. This is another proof of Theorem \ref{prop:p-kallenberg}, since (\ref{eq:pcid})  implies the partially c.i.d. condition.

\begin{exa}\label{ex:gauss}\rm
{\bf [A Gaussian partially c.i.d. model]}.
Finding the explicit expression of the limit laws $\alpha_i$ is, in general, a difficult task. Yet, we present here a partially c.i.d. system for which the limit law is available and has a parametric expression. Details on computations are given in the Appendix. 

Think of unevenly spaced, synchronous observations on  $K$ variables. Let $T_1<T_2<\dots$ denote the times at which observations are taken. Assume that $T_n<\infty$, $\mathbb P$-a.s. for every $n\geq 1$. For every $n\geq 0$, let $t_n=T_{n+1}-T_{n}$ denote the inter-arrival times ($T_0=0$), and let $X_{n,i}$ denote the observation on the $i$-th variable at time $T_n$, $i\in I=\{1,\dots,K\}$. Assume that the $X_{n,i}$ ($i\in I$) are conditionally independent given $\mathcal G_{n-1}=\sigma(T_{1:n},\bX_{1:n-1}^I)$, and that $T_{n+1}$ is conditionally independent of $\bX_n^I:=(X_{n,i}:i=1,\dots,K)$, given  $\mathcal G_{n-1}$, with
\begin{equation} \label{eq:gaussPred}
\left\{
\begin{aligned}
&X_{1,j}\sim Q_{0,j}=\Norm(\mu_{1,j},\sigma^2_{1,j}), \\
& X_{n,j}\mid\mathcal G_{n-1} \sim Q_{n-1,j}=\Norm(\mu_{n,j},\sigma^2_{n,j}),\quad n\geq 2
\end{aligned}\right.
\end{equation}
where $\mu_{1,j}\in\mathbb R$, $\sigma^2_{1,j}>0$ and $\mu_{n,j}$ is given by  the \textit{last tick} (or {\it piecewise-constant}) {\it interpolation scheme} of the process (see e.g. \cite{genccay2001}, \cite{hayashi2005}):
\[
\mu_{n,j}=\frac{t_0\mu_{1,j}+\sum_{k=1}^{n-1}t_kX_{k,j}}{T_{n}},
\quad n\geq 2.
\]
A sufficient condition for $[X_{n,i}]$ to be partially $\mathcal G$-c.i.d. is
\[
\sigma^2_{n,j}=\prod_{k=1}^{n-1}\left[1-t_k^2/T_{k+1}^2\right]\sigma^2_{1,j} , \quad n\geq 2.  
\]
Then, defining $\lambda_n=t_n/T_{n+1}$, we can write
\begin{equation} \label{eq:gaussrecurs}
\left\{
\begin{aligned}
&\mu_{n,j}=(1-\lambda_{n-1})\mu_{n-1,j}+\lambda_{n-1}X_{n-1,j}\\
&\sigma^2_{n,j}=(1-\lambda_{n-1}^2)\sigma^2_{n-1,j}
\end{aligned}\right. \quad n\geq 2. 
\end{equation}
By Proposition \ref{prop:strong}, there exist $K$ random probability measures $\alpha_1,\dots, \alpha_K$ such that $(Q_{n,1},\dots, Q_{n,K})\rightarrow (\alpha_1,\dots,\alpha_K)$ $\mathbb P$-a.s.,   and 
$ \mathbb P[\bX_{n}^I\in\cdot \mid
\bX_{1:n-1}^I, T_{1:n}]\rightarrow \alpha_1\times\dots\times\alpha_K$. 
It can be proved that
\[
\alpha_i=\Norm(\mu_i,\gamma \sigma^2_{1,i}) \quad (i=1,\dots,K),
\]
with $\gamma=\prod_{k=1}^\infty (1-\lambda_k^2)$ and
$$ 
\mu_i \mid \gamma \stackrel{ind}{\sim}\Norm(\mu_{1,i},(1-\gamma)\sigma_{1,i}^2). 
$$
Hence, the random probability measures $\alpha_1,\dots,\alpha_K$ are not independent, unless $\gamma$ has a degenerate distribution. For example, they are not independent when the observations arrive at a Poisson process rate. 
\end{exa}

\section{A strong law of large numbers} \label{sec:SLLN}

In this section, we provide a strong law of large numbers (SLLN) for $f(\bX_n^I)$, with $f$ measurable and  $[X_{n,i}]$ a partially c.i.d. array with directing random measures $(\alpha_i,i\in I)$. Notice that the sequence  $(f(\bX_n^I))$ is not c.i.d. in general. 
Since the SLLN is based on the convergence of $\mathbb E[f(\bX_{n+1}^I)\mid \mathcal G_n]$, we first give sufficient conditions for
\begin{equation}
\label{eq:conve_cond_expe}
\mathbb E[f(\bX_{n+1}^I)\mid \mathcal G_n]\rightarrow \int f d\alpha  \quad \mathbb P\mbox{-a.s.},
\end{equation}
where $\alpha=\mathlarger{\mathlarger{\mathlarger{\times}}}_{i \in I} \alpha_i$. We know, from Proposition \ref{prop:strong}, that  (\ref{eq:conve_cond_expe}) holds if $f$ is bounded and measurable. The following proposition
shows that a necessary and sufficient condition for $(f(\bX_{n}^I))$ to satisfy (\ref{eq:conve_cond_expe})
is predictive uniform integrability.
We say that a sequence $(Y_n)$ of real-valued random variables is $\mathcal G$-{\em predictive uniformly integrable} if $Y_n$ is integrable for every $n\geq 1$ and
	\begin{equation}
	\label{eq:unif_integr2}
	\sup_n	\mathbb E[|Y_{n+1}|\mathbbm{1}_{\{|Y_{n+1}|>k\}}\mid\mathcal G_n]\rightarrow 0\quad \mathbb P\mbox{-a.s.}\quad \mbox{ as }k\rightarrow\infty.
	\end{equation}
	
\begin{prop}
	\label{th:convpred}
 Let  $(\bX_n^I)$  be a partially $\mathcal G$-c.i.d. process  with directing random measures $(\alpha_i, i \in I)$ and let $f:\mathbb X^I\rightarrow \mathbb R$ be a measurable function. Then,  the following conditions are equivalent
 \begin{itemize}
 	\item[(i)] $f(\bX_n^I)$ is $\mathcal G$-predictive uniformly integrable;
 	\item[(ii)] $\mathbb E[f(\bX_{n+1}^I)\mid \mathcal G_n] \rightarrow \int f d\alpha<\infty$,
 	$\mathbb P$-a.s. as $n\rightarrow\infty$, where $\alpha=\mathlarger{\mathlarger{\mathlarger{\times}}}_{i \in I} \alpha_i$.
	\end{itemize}
\end{prop}

\noindent \textsc{Proof}. Let $Y_n=f(\bX_n^I)$ and, for every fixed $k\in\mathbb N$, let $f_k=f\mathbbm{1}_{\{f\leq k\}}$.\\
(i)$\Rightarrow$(ii). We prove the result for nonnegative $f$. The general case can be obtained by linearity.
 We know, from Proposition \ref{prop:strong}, that $\mathbb E[Y_{n+1}\mathbbm{1}_{\{Y_{n+1}\leq k\}}\mid\mathcal G_n]$ converges $\mathbb P$-a.s. to $\int f_kd\alpha$. Hence,
\[
\liminf_n\mathbb E[Y_{n+1}\mid\mathcal G_n]\geq \liminf_n\mathbb E[Y_{n+1}\mathbbm{1}_{\{Y_{n+1}\leq k\}}\mid\mathcal G_n]\geq\int f_kd\alpha\quad \mathbb P\mbox{-a.s.}
\]
Since $\int f_kd\alpha \rightarrow\int fd\alpha$ $\mathbb P$-a.s as $k\rightarrow\infty$, $$\liminf_n\mathbb E[Y_{n+1}\mid\mathcal G_n]\geq \int fd\alpha\quad \mathbb P\mbox{-a.s.}$$
To show the reverse inequality, let $N_1\in\mathcal F$ be such that $\mathbb P[N_1]=0$ and $\mathbb E[Y_{n+1}\mathbbm{1}_{\{Y_{n+1}\leq k\}  }\mid\mathcal G_n](\omega)$ converges to $\int f_kd\alpha(\omega)$ for every $k\in \mathbb N$ and $\omega\in N_1^c$. Let $N_2$ be a null set such that the convergence in (\ref{eq:unif_integr2}) holds for every $\omega\in N_2^c$. Then, for every $\omega\in N_2^c$ and every $\epsilon>0$, there exists $k\in\mathbb N$ such that, 
$\sup_n \mathbb E[Y_{n+1}\mathbbm{1}_{\{Y_{n+1}>k  \}  }\mid\mathcal G_n](\omega)<\epsilon$. Hence, for every $\omega\in N_1^c\cap N_2^c$,
\[
\limsup_n \mathbb E[Y_{n+1}\mid \mathcal G_n](\omega)\leq \limsup_n \mathbb E[Y_{n+1}\mathbbm{1}_{\{Y_{n+1}\leq k  \}  }\mid \mathcal G_n](\omega)+\epsilon\leq \int f_kd\alpha(\omega)+\epsilon\leq \int fd\alpha(\omega)+\epsilon.\]
Since this holds for every $\epsilon$, $\limsup_n \mathbb E[Y_{n+1}\mid \mathcal G_n]\leq \int fd\alpha$ $\mathbb P$-a.s.
Hence, $ E[Y_{n+1}\mid\mathcal G_n]$ converges $\mathbb P$-a.s. to $\int fd\alpha$. To show that$\int fd\alpha<\infty$ $\mathbb P$-a.s., notice that, for every $\omega\in N_1^c\cap N_2^c$,
$$
\limsup_n\mathbb E[Y_{n+1}\mid\mathcal G_n](\omega)\leq \limsup_n\mathbb E[Y_{n+1}\mathbbm{1}_{\{Y_{n+1}\leq k\}}\mid\mathcal G_n](\omega)+\sup_n\mathbb E[Y_{n+1}\mathbbm{1}_{\{Y_{n+1}> k\}}\mid\mathcal G_n](\omega)\leq k+\epsilon.
$$
(ii)$\Rightarrow$(i). Since $\mathbb E[| Y_{n+1}| \mathbbm{1}_{\{| Y_{n+1}| >k\}}\mid \mathcal G_n]=\mathbb E[Y^+_{n+1} \mathbbm{1}_{\{Y^+_{n+1}>k\}}\mid \mathcal G_n]+\mathbb E[Y^-_{n+1}\mathbbm{1}_{\{Y^-_{n+1}>k\}}\mid \mathcal G_n]$, it is sufficient to prove predictive uniform integrability for non negative $f$.
For every $k$ and $n$,
\[
\mathbb E[Y_{n+1}\mathbbm{1}_{\{Y_{n+1}>k \}}\mid \mathcal G_n]\leq \left| \mathbb E[Y_{n+1}\mid\mathcal G_n]-\int fd\alpha\right|+\left| \int fd\alpha-\int f_k d\alpha\right| +\left| \int f_k d\alpha -\mathbb E[Y_{n+1}\mathbbm{1}_{\{Y_{n+1}\leq k \}}\mid\mathcal G_n]\right|.
\]
Let $N$ be such that $\mathbb P[N]=0$ and, for every $\omega\in N^c$, $\mathbb E[Y_{n+1}\mid \mathcal G_n]\rightarrow \int fd\alpha$ 
and $\mathbb E[Y_{n+1}\mathbbm{1}_{\{Y_{n+1}\leq k\}}\rightarrow \int f_k d\alpha$ for every $k$. Since $\int fd\alpha<\infty$, $\mathbb P$-a.s., then, for every fixed $\epsilon$ and $\omega\in N^c$, there exists $k$ such that $\mid \int f d\alpha(\omega)-\int f_kd\alpha(\omega) \mid <\epsilon/3$. For such a $k$, let $n_0$ be such that, for every $n\geq n_0$,
$|\mathbb E[Y_{n+1}\mid \mathcal G_n](\omega)-\int fd\alpha(\omega)| <\epsilon/3$ and $|\mathbb E[Y_{n+1}\mathbbm{1}_{\{Y_{n+1}\leq k\}}\mid\mathcal G_n](\omega)-\int f_kd\alpha(\omega)| <\epsilon/3$. Then,
\[
\sup_{n\geq n_0}\mathbb E[Y_{n+1}\mathbbm{1}_{\{Y_{n+1}>k \}}\mid\mathcal G_n](\omega)<\epsilon.
\]
\noindent $\square$

By Proposition \ref{th:convpred}, predictive uniform integrability implies (\ref{eq:conve_cond_expe}). If $f(\bX_n^I)$ is also dominated in $L^1$, then we get the SLLN.
We say that a sequence $(Y_n)$ of real-valued random variables is \textit{dominated in} $L^1$ if
	\begin{equation}
	\label{eq:unifint1}
	\int_0^\infty \sup_n\mathbb P[|Y_n|>x]dx<\infty.
	\end{equation}
$(Y_n)$ is dominated in $L^1$ if and only if $|Y_n|$ is stochastically dominated by an integrable random variable $Y$, for any $n$.
\begin{thm} {\em (SLLN)}
	\label{th:slln}
	 Let  $(\bX_n^I)$  be a partially $\mathcal G$-c.i.d. process  with directing random measures $(\alpha_i, i \in I)$, and let  $f:\mathbb X^I\rightarrow \mathbb R$ be a measurable function. If $(f(\bX_n^I))$  is $\mathcal G$-predictive  uniformly integrable and dominated in $L^1$, then, as $n\rightarrow\infty$,
	 \[
	\frac{1}{n}\sum_{k=1}^n f(\bX_k^I)\rightarrow \int fd\alpha , \quad \mathbb P\mbox{-a.s.},
	\]
where $\alpha= \mathlarger{\mathlarger{\mathlarger{\times}}}_{i\in I}\alpha_i$.
\end{thm}

\noindent \textsc{Proof}.  We prove the theorem for nonnegative $f$. The general case can be obtained by linearity.\\
	For every $n\geq 1$, let $Y_n=f(\bX_n^I)$ and  $Y_n^*=Y_n \mathbbm{1}_{\{Y_n\leq n\}}$. Furthermore, let $F(x)=\inf_n \mathbb P[ Y_{n} \leq x]$. Then,
		\[
		\sum_{n=1}^\infty \mathbb P[Y_n\neq Y_n^*]
		=	\sum_{n=1}^\infty \mathbb P[Y_n>n]\\
		\leq \sum_{n=1}^\infty (1-F(n)).
		\]
		By (\ref{eq:unifint1}), $ \sum_{n=1}^\infty (1-F(n))<\infty$. Hence, $\mathbb P[Y_n\neq Y_n^*\;i.o.]=0$. It is therefore sufficient to prove that
		$\frac{1}{n}\sum_{k=1}^n Y_k^*\rightarrow \int fd\alpha$, $\mathbb P$-a.s. Let
	\[
	M_n=\sum_{k=0}^{n-1}\frac{Y_{k+1}^*-\mathbb E[Y_{k+1}^*\mid\mathcal G_k]}{k+1}.
	\]
	Then $M_n$ is a $\mathcal G$-martingale and, by (\ref{eq:unifint1}),
	\[
	\begin{aligned}
	\mathbb E[M_n^2]&=\sum_{k=0}^{n-1}\frac{\mathbb E[Y_{k+1}^*-\mathbb E[Y_{k+1}^*\mid\mathcal G_k]]^2}{(k+1)^2}\leq \sum_{k=1}^{\infty}\frac{\mathbb E[Y_{k}^{*2}]}{k^2}\leq 2\sum_{k=1}^{\infty}\frac{\int_0^{k}x\mathbb P[Y_{k}>x]dx}{k^2}\\&	\leq 2\sum_{k=1}^\infty \frac{\sum_{h=1}^k h (1-F(h-1))}{k^2}
 \leq 2\sum_{h=1}^\infty h(1-F(h-1))\sum_{k=h}^\infty\frac{1}{k^2}
\leq 4\sum_{h=1}^\infty (1-F(h-1)),
	\end{aligned}
	\]
	which is finite, by (\ref{eq:unifint1}).
  It follows that $M_n$ converges $\mathbb P$-a.s. Hence,
	\[
	\sum_{k=0}^{\infty}\frac{Y_{k+1}^*-\mathbb E[Y_{k+1}^*\mid\mathcal G_k]}{k+1}<\infty \, \quad  \mathbb P\mbox{-a.s.}
	\]
	By Kronecker's Lemma,
	\[
	\frac{1}{n}\sum_{k=0}^{n-1}\left(Y_{k+1}^*-\mathbb E[Y_{k+1}^*\mid\mathcal G_k]\right)\rightarrow 0 \,   \quad \mathbb P\mbox{-a.s.}
	\]
	as $n\rightarrow\infty$. The proof is complete if we can show that 	
	\begin{equation}\label{eq:conv_emp_pred}
	\mathbb E[Y_{n+1}^*\mid \mathcal G_n]\rightarrow \int fd\alpha \,  \quad  \mathbb P\mbox{-a.s.}
	\end{equation}
Since $(Y_n)$ is $\mathcal G$-predictive uniformly integrable, by Proposition \ref{th:convpred}  $\mathbb E[Y_{n+1}\mid\mathcal G_n]$ converges to $\int fd\alpha$, $\mathbb P$-a.s. On the other hand, $\mathcal G$-predictive uniform integrability implies that $\mathbb E[Y_{n+1}\mathbbm{1}_{\{Y_{n+1}>n+1\}}\mid\mathcal G_n]$ converges to zero $\mathbb P$-a.s. Hence, (\ref{eq:conv_emp_pred}) holds.\\
\noindent $\square$
	
The following propositions give sufficient conditions for a sequence $(Y_n)$ to be  $\mathcal G$-predictive uniformly integrable and dominated in $L^1$.
\begin{prop}
	\label{prop:suffcond_dom}
	Let $(Y_n)$ and $(Y'_n)$ be sequences of real-valued random variables,
		such that $|Y_n|\leq  |Y'_n|$, $\mathbb P$-a.s. for every $n$.
		\begin{itemize}
\item[(i)]	If $(Y'_n)$ is $\mathcal G$-predictive uniformly integrable, then $(Y_n)$ is $\mathcal G$-predictive uniformly integrable. 
\item[(ii)] If $(Y'_n)$ is dominated in $L^1$, then $(Y_n)$ is dominated in $L^1$.
\end{itemize}
\end{prop}

\begin{prop}
	\label{prop:suffcond_r}
	Let $(Y_n)$ be a sequence of real-valued random variables, and  $\mathcal G=(\mathcal G_n)$ a filtration. 
	\begin{itemize}
		\item[(i)] If $
	\sup_n	\mathbb E[|Y_{n+1}|^r\mid\mathcal G_n]<\infty$
$\mathbb P$-a.s. for some $r>1$,
	then  $(Y_n)$ is $\mathcal G$-predictive uniformly integrable.
\item[(ii)] If
  $\sup_n\mathbb E[|Y_{n}|^r]<\infty$,
	for some $r>1$, then $(Y_n)$ is dominated in $L^1$.
			\end{itemize}
		\end{prop}
The proofs of Propositions \ref{prop:suffcond_dom} and \ref{prop:suffcond_r} follow from standard arguments.

\begin{prop}
	\label{prop:suffcond_cid}
	If $(Y_n)$ is a $\mathcal G$-c.i.d. sequence of integrable real-valued random variables, then $(Y_n)$ is $\mathcal G$-predictive uniformly integrable and dominated in $L^1$.
\end{prop}
\noindent \textsc{Proof}. It is enough to prove the thesis for non-negative $Y_n$. 
The result can be extended by linearity.
Since $\int_0^\infty \sup_n \mathbb P[Y_n >x]dx=\mathbb E[Y_1]<\infty$, $(Y_n)$ is dominated in $L^1$.
To prove predictive uniform integrability, we use Proposition \ref{th:convpred} with $\bX_n^I=Y_n$ and $f$ the identity function, and show that
$\mathbb E[Y_{n+1}\mid \mathcal G_n]$ converges almost surely to $\int y \alpha(dy)$, where $\alpha $ is the directing random measure of $(Y_n)$.
We first show that  $(\mathbb E[Y_{n+1}\mid\mathcal G_n])_{n\geq 0}$ is a uniformly integrable martingale. The martingale property follows from Proposition \ref{prop:kallenberg2} if the distribution of $Y_1$ has bounded support, and can be extended to the general case by truncation and the monotone convergence theorem. Uniform integrability follows from $\mathbb E[\mathbb E[Y_{n+1}\mid \mathcal G_n]\mathbbm 1_{\{\mathbb E[Y_{n+1}\mid\mathcal G_n]>a\}}]\leq \mathbb E[Y_{n+1}\mathbbm 1_{\{Y_{n+1}>a\}}]$ $\mathbb P$-a.s., which holds by conditional Jensen's inequality, and from uniform integrability of $(Y_n)$. Being a uniformly integrable martingale,  $(\mathbb E[Y_{n+1}\mid\mathcal G_n])_{n\geq 0}$  converges $\mathbb P$-a.s. to an integrable random variable $V$ such that $\mathbb E[V\mid \mathcal G_n]=\mathbb E[Y_{n+1}\mid\mathcal G_n]$, $\mathbb P$-a.s.
To prove that $V=\int y \alpha(dy)$, let us define, for every  $k\in\mathbb N$, $b_k(y)=y \mathbbm{1}_{\{y\leq k\}}$. Then,  $\mathbb E[b_k(Y_{n+1})\mid\mathcal G_n]$ converges $\mathbb P$-a.s. to $V_{k}=\int b_k(y)\alpha(dy)$ by  Proposition \ref{prop:strong}, and $\mathbb E[V_k\mid \mathcal G_n]=\mathbb E[b_k(Y_{n+1})\mid\mathcal G_n]$ $\mathbb P$-a.s.
Let us prove that $V_{k}$ converges to $V$ $\mathbb P$-a.s. as $k\rightarrow\infty$.
 Indeed, $\mathbb E[V-V_{k}]=\mathbb E[Y_1\mathbbm{1}_{\{Y_1> k\}}]\rightarrow 0$, as $k\rightarrow\infty$, since $Y_1\mathbbm{1}_{\{Y_1> k\}}\rightarrow 0$ $\mathbb P$-a.s. and is bounded by $Y_1$, which is integrable. Thus, the sequence of non negative random variables $(V-V_{k})$ converges to zero in probability, as $k\rightarrow\infty$. Being a monotone sequence, $(V-V_{k})$ converges $\mathbb P$-a.s. to zero, as $k\rightarrow\infty$.
   On the other hand, $V_{k}=\int b_k(y)\alpha(dy)\rightarrow \int y\alpha(dy)$, $\mathbb P$-a.s. as $k\rightarrow \infty$, where the last convergence follows from the monotone convergence theorem. Hence, $V=\int y\alpha(dy)$ $\mathbb P$-a.s. \\
\noindent $\square$

The results are illustrated in the following examples.
  \begin{exa} \rm
  	\label{exa:1}
Let $[X_{n,i}]$ be a partially c.i.d. array, with $i=1,\dots,p$, such that  $\mathbb E[|X_{1,i}|^p]<\infty$
 for all $i$.  Then,
 \begin{equation}
 \label{eq:example1_ssln}
 \frac{1}{n}\sum_{k=1}^n \prod_{i=1}^pX_{k,i}\rightarrow \prod_{i=1}^p\int x\alpha_i(dx)\quad \mathbb P\mbox{-a.s.},
 \end{equation}
 where $(\alpha_1,\dots,\alpha_p)$ are the directing random measures of $[X_{n,i}]$.
To prove it, we can write the right-hand side of (\ref{eq:example1_ssln}) as $\int f(x_1,\dots,x_p)\alpha(dx_1,\dots,dx_p)$, with $f(x_1,\dots,x_p)=\prod_{i=1}^px_i$ and $\alpha=
\Times_i\alpha_i$ and apply the strong law of large numbers to  $Y_n=f(X_{n,1},\dots,X_{n,p})$. Notice that $|Y_n|\leq \sum_{i=1}^p|X_{n,i}|^p/p$ with $|X_{n,i}|^p$ integrable and c.i.d. with respect to the natural filtration of $[X_{n,i}]$. By Propositions \ref{prop:suffcond_dom} and \ref{prop:suffcond_cid}, $(Y_n)$ is dominated in $L^1$ and predictive uniformly integrable with respect to the natural filtration of $[X_{n,i}]$. By Theorem \ref{th:slln}, $\sum_{j=1}^nY_j/n$ converges $\mathbb P$-a.s. to $\int fd\alpha$. \\
\noindent $\square$
\end{exa}

\medskip

\begin{exa} \rm
	\label{exa:2}
 Let $[X_{n,i}]$ be a partially c.i.d. array of positive  random variables, with $i=1,\dots,p$, such that
 $\mathbb E[|\log X_{1,i}|]<\infty$
for all $i$. Then 
 \begin{equation*}
 \frac{1}{n}\sum_{k=1}^n \log(\sum_{i=1}^pX_{k,i})\rightarrow \int \log(\sum_{i=1}^px_i)\alpha_1(dx_1)\dots\alpha_p( dx_p)\quad \mathbb P\mbox{-a.s.},
 \end{equation*}
 where $(\alpha_1,\dots,\alpha_p)$ are the directing random measures of $[X_{n,i}]$. To prove this, let $Y_n=f(X_{n})$, with $f(x_1,\dots,x_p)=\log(\sum_{i=1}^px_i)$. Notice that
 \[
\log(\sum_{i=1}^p X_{n,i})\leq \left\{
 \begin{array}{ll}
 \max_{1\leq i\leq p}|\log X_{n,i}|&\mbox{ if }\sum_{i=1}^pX_{n,i}<1\\
 \sum_{1\leq i\leq p}\log (X_{n,i}+2)\mathbbm{1}_{\{X_{n,i}<2\}} +
 \sum_{1\leq i\leq p}\log X_{n,i}\mathbbm{1}_{\{X_{n,i}\geq 2\}}&\mbox{ if }\sum_{i=1}^pX_{n,i}\geq 1.
  \end{array}
 \right.
 \]
Therefore,  $|Y_n|\leq \sum_{i=1}^p|\log X_{n,i}|+p\log 4$. Since $\log X_{n,i}$ is integrable and c.i.d. with respect to the natural filtration of $[X_{n,i}]$, by Propositions \ref{prop:suffcond_dom} and \ref{prop:suffcond_cid}, $(\log(\sum_{i=1}^p X_{n,i}))$ is predictive uniformly integrable with respect to the natural filtration of $[X_{n,i}]$ and dominated in $L^1$. By theorem \ref{th:slln}, $\sum_{k=1}Y_k/n$ converges $\mathbb P$-a.s. to $\int fd\alpha$.\\ \noindent $\square$
 \end{exa}

\section{Central limit theorems}   \label{sec:CLT}
This section deals with the central limit problem for
	$(f_i(X_{n,i}),i\in I)$, when $[X_{n,i}]$ is partially c.i.d. and the $f_i$ are real valued functions defined on $\mathbb X$. Notice that, if $[X_{n,i}]$ is partially $\mathcal G$-c.i.d., the array $[f_i(X_{n,i})]$ is also partially $\calG$-c.i.d. Therefore, we can rename $f_i(X_{n,i})$ as $X_{n,i}$, to simplify the notation. Thus, throughout the current section, the $X_{n,i}$ are assumed to be \textit{real-valued} random variables.

We give two Central Limit Theorems (CLTs) for the empirical sums
	$(\sum_{k=1}^n X_{k,i},  i\in I)$, suitably centered and scaled.
	The choice of the centering focuses on the agreement between predictions and empirical means.
Namely, we prove CLTs for the scaled cumulative forecast errors $\sum_{k=1}^n (X_{k,i} - \mathbb E[X_{k,i} \mid \calG_{k-1}])$ and for the  approximation errors  $\sum_{k=1}^n X_{k,i}/n - \E[X_{n+1} \mid \calG_n]$. In statistical inference, these forms of CLTs may provide a basis for model checking and for the approximation of predictions, in problems such as Bayesian forecasting with large samples, discrete time filtering and sequential procedures, where prediction is the main focus of interest but exact computations are costly. It is worth noticing that we deal with the {\em joint} scaled differences, for $i \in I$. Thus, the CLTs will involve a multivariate Gaussian distribution, when $I$ is finite, or a Gaussian measure on ${\mathbb R}^I$,  when $I$ is countable.

As in the case of exchangeable and c.i.d. sequences, CLTs for partially c.i.d. sequences are given in terms of stable convergence.  Stable convergence is stronger than convergence in distribution and weaker than convergence in probability.  We recall some basic notions regarding stable convergence, and refer the reader to \cite{Aldous85} 
and \cite{hausler2015} for further reading. Roughly speaking, stable convergence means convergence of the conditional distributions. A sequence $(Z_n)$ of random variables  defined on $(\Omega, \calF, \bP)$ and  taking values in a Polish space $\mathbb X$ {\em converges stably} if, for each non-null event $A$, the conditional distribution of $Z_n$ given $A$ converges weakly. If $Z_n$ converges stably, then there exists a {\em representing} random measure  $\beta(\cdot; \omega)$ such that one can write
$$
 \bP[Z_n \in \cdot\mid A] \,\bP[A] \rightarrow \int_A \beta(\cdot ; \omega) \, \bP(d \omega), $$
and $Z_n$ is said to {\em converge stably with representing measure}  $\beta(\cdot; \omega)$.
For $A=\Omega$, this implies that the limit measure is a mixture, 
$$ \bP[Z_n \in \cdot] \rightarrow \int_{\Omega} \beta(\cdot; \omega) \, \bP(d \omega).
$$
For example, if $\beta(\cdot; \omega)$ is a Gaussian probability law $\Norm(\cdot; 0, \sigma^2(\omega))$, with zero mean and random variance $\sigma^2(\omega)>0$, one has, for every $B$,
$$
P[Z_n \in B] \rightarrow \int_\Omega \Norm(B; 0, \sigma^2(\omega)) \, \bP(d \omega); $$
moreover, for any $A$ and $B$,
$$ \bP[Z_n \in B, \sigma^2 \in A] \rightarrow \int_{\{\sigma^2(\omega) \in A\}} \Norm(B ; 0, \sigma^2(\omega)) \bP(d \omega).$$
If $Z_n$ converges stably, one can extend the space $(\Omega, {\cal F}, \bP)$  to construct a \lq \lq limit'' variable $Z$ such that the representing measure is a regular conditional distribution for $Z$ given ${\cal F}$ (see \cite{Aldous85}, page 57). Then another relevant implication is that, if $Z_n \rightarrow Z$ stably,
and $Y_n \rightarrow Y$ in probability, then $(Z_n, Y_n)$ converges in distribution to $(Z, Y)$. This result allows a generalized version of the  Cram\'{e}r-Slutzky theorem, covering the case when $Y$ is a random variable; see \cite{hausler2015}.

 The results in this section extend ideas in \cite{BertiPratelliRigo2004LimThForIID} and are based on the following central limit theorem for martingale differences. 

\begin{thm}{\em (\cite{HH80Martingale}, Theorem 3.2)}\label{lemma:clt}
	Let $[Z_{n,k}]$, $k=1,\dots,k_n,n\geqslant 1$ be an array of real, square-integrable random variables, with $k_n\rightarrow \infty$, and let $S_n=\sum_{k=1}^{k_n}Z_{n,k}$. For all $n$, let $\mathcal{F}_{n,0}\subset \mathcal{F}_{n,1}\subset\dots\subset\mathcal{F}_{n,k_n}\subset\mathcal{F}$ be $\sigma$-fields with $\mathcal{F}_{n,0}=\{\emptyset,\Omega\}$.   If the following conditions hold:
\begin{equation}
\label{enum:clt_1}
\sigma(Z_{n,k})\subset \mathcal{F}_{n,k}\subset \mathcal{F}_{n+1,k};
\end{equation}
\begin{equation}
\label{enum:clt_2}
 \mathbb{E}[Z_{n,k}\mid\mathcal{F}_{n,k-1}]=0;
 \end{equation}
\begin{equation}
\label{enum:clt_3}
\max_{1\leqslant k \leqslant k_n}|Z_{n,k}|\stackrel{P}{\rightarrow} 0 ;
\end{equation}
\begin{equation}
\label{enum:clt_4}
\underset{n}{\sup}\;\;\mathbb{E}[\max_{1\leqslant k\leqslant k_n}Z_{n,k}^2]<\infty;
\end{equation}
\begin{equation}
\label{enum:clt_5}
\sum_{k=1}^{k_n}Z_{n,k}^2\stackrel{P}{\rightarrow} L,
\end{equation}
then	$(S_n)$ converges stably with representing measure  $\Norm(0,L)$.
\end{thm}

The first limit result that we consider deals with the sums of the \textit{forecast errors} 
		\begin{equation}
		\label{eq:Uki}
		U_{n,i}=X_{n,i}-\mathbb E[X_{n,i}\mid \mathcal G_{n-1}], \quad \quad  n \geq 1, i \in I. 
		\end{equation}	
The limit distribution turns out to be a mixture of $\Norm(0,\Sigma)$ probability laws, where $\Norm(0,\Sigma)$ denotes the centered Gaussian measure on $\mathbb R^I$ with correlation operator $\Sigma$ (see \cite{daletskii1991measures}).  We  denote by $(e_i, i\in I)$ the canonical basis of $\mathbb R^I$.
\begin{thm}\label{thm:clt1}
		Let $[X_{n,i}]$ be a partially $\calG$-c.i.d. array of real-valued and square-integrable random variables, with directing random measures $(\alpha_i:i\in I)$. Let 
		\[
		S_{n,i}=\frac{\sum_{k=1}^n\left(X_{k,i}-\mathbb E[X_{k,i}\mid \mathcal G_{k-1}]\right)}{\sqrt{n}}\quad\quad\quad n\geq 1,i\in I.
		\]
		Then, the sequence $\left(\mathbf S_{n}^I\right)$
	 converges stably  with representing measure {\rm $\Norm(0,\Sigma)$}, where $\Sigma$ is a diagonal operator and $\Sigma(e_i,e_i)=\sigma_{\alpha_i}^2$, with $\sigma^2_{\alpha_i}=\int x^2\alpha_i(dx)-(\int x \alpha_i(dx))^2$. 	 
\end{thm}
	
\noindent \textsc{Proof.} It is sufficient to show that, for every finite $J\subset I$ and for every  $(a_i,i\in J)\in \mathbb{R}^J$, $\sum_{i\in J}a_i S_{n,i}$ converges stably, with representing measure $\Norm(0,\sum_{i,j\in J}a_{i}a_j\Sigma_{i,j})$, where $\Sigma_{i,j}=\Sigma(e_i,e_j)$ (see \cite{hausler2015},  Corollary 3.19 and Proposition 3.22, and \cite{daletskii1991measures}). 
We use Theorem \ref{lemma:clt}, letting $k_n=n$, $\mathcal{F}_{n,k}=\mathcal{G}_k$ and $Z_{n,k}=n^{-1/2}\sum_{i\in J}a_iU_{k,i}$, with $(U_{k,i})$ as in (\ref{eq:Uki}). Then $\sum_{i\in J}a_i S_{n,i}=\sum_{k=1}^n Z_{n,k}$ and checking \eqref{enum:clt_1}-\eqref{enum:clt_5} is sufficient to obtain the thesis.
	Conditions \eqref{enum:clt_1} and \eqref{enum:clt_2} are immediate.
	To prove  \eqref{enum:clt_3}, notice that 
	$$
	\max_{1\leqslant k \leqslant k_n}|Z_{n,k}|\leq \sum_{i\in J}|a_i|\max_{1\leqslant k \leqslant n}|U_{k,i}|/\sqrt n
	\leq \sum_{i\in J}|a_i|\max_{1\leqslant k \leqslant n}|X_{k,i}|/\sqrt n +
	\sum_{i\in J}|a_i|\max_{1\leqslant k \leqslant n}\mathbb E[|X_{k,i}|\mid \mathcal G_{k-1}]/\sqrt n.
	$$ 
	Furthermore,
	$$
	\mathbb P[\max_{1\leqslant k \leqslant n}|X_{k,i}|>\epsilon \sqrt n]\leq \sum_{k=1}^n \mathbb P[|X_{k,i}|>\epsilon \sqrt n]\leq \frac{1}{\epsilon^2}\mathbb E[X_{1,i}^2\mathbbm{1}_{\{X_{1,i}^2>\epsilon^2 n\}}]
	$$
	$$
	\mathbb P[\max_{1\leqslant k \leqslant n}\mathbb E[|X_{k,i}|\mid\mathcal G_{k-1}]>\epsilon \sqrt n]\leq \sum_{k=1}^n  \mathbb P[\mathbb E[|X_{k,i}|\mid\mathcal G_{k-1}]>\epsilon \sqrt n]\leq \frac{1}{\epsilon^2}\mathbb E[X_{1,i}^2\mathbbm{1}_{\{X_{1,i}^2>\epsilon^2 n\}}],
	$$
	which converges to zero, as $n \rightarrow \infty$. 
	To prove  \eqref{enum:clt_4}, notice that
	$$
	\sup_n\;\;\mathbb E[\max_{1\leq k\leq n}Z_{n,k}^2]\leq \sup_n\frac{1}{n}\sum_{i,j\in J}|a_ia_j|\mathbb E[\max_{1\leq k\leq n}U_{k,i}^2]
	$$
	and
	$$
	\sup_n\frac{1}{n}\mathbb E[\max_{1\leq k\leq n}U_{k,i}^2]\leq \sup_n\frac{1}{n}\sum_{k=1}^n \mathbb E[(X_{k,i}-\mathbb E[X_{k,i}\mid\mathcal G_{k-1}])^2]\leq \sup_n\frac{1}{n}\sum_{k=1}^n \mathbb E[X_{k,i}^2]<\infty.
	$$
	To prove \eqref{enum:clt_5}, we write
	\begin{equation}
	\label{eq:Z^2}
	\sum_{k=1}^{n}Z_{n,k}^2=\frac{1}{n}\sum_{k=1}^{n}\left(\sum_{i\in J}a_i U_{k,i}\right)^2=\sum_{i,j\in J}a_ia_j\frac{1}{n}\sum_{k=1}^{n} U_{k,i} U_{k,j}
	\end{equation}
	and
	$$
	\frac{1}{n}\sum_{k=1}^{n} U_{k,i} U_{k,j}=\frac{1}{n}\sum_{k=1}^{n} X_{k,i}X_{k,j}-
	\frac{1}{n}\sum_{k=1}^{n} \mathbb E[X_{k,i}\mid\mathcal G_{k-1}]\mathbb E[X_{k,j}\mid\mathcal G_{k-1}]+E_{n,i,j}+E_{n,j,i},
	$$
	where, for every $i$ and $j$ in $J$,
	$E_{n,i,j}=\frac{1}{n}\sum_{k=1}^{n} (X_{k,i}-\mathbb E[X_{k,i}\mid\mathcal G_{k-1}])\mathbb E[X_{k,j}\mid\mathcal G_{k-1}]$. Since  $\frac{1}{n}\sum_{k=1}^{n} (X_{k,i}-\mathbb E[X_{k,i}\mid\mathcal G_{k-1}])$ converges $\mathbb P$-a.s. to zero, and $\mathbb E[X_{k,j}\mid\mathcal G_{k-1}]$ converges $\mathbb P$-a.s. to $\int x\alpha_i(dx)$, then $E_{n,i,j}$ converges to zero, $\mathbb P$-a.s. for every $i,j\in J$.
	Furthermore, by Theorem \ref{th:slln} and the computations in Example \ref{exa:1}, 	
	 $\frac{1}{n}\sum_{k=1}^{n} X_{k,i}X_{k,j}$ converges $\mathbb P$-a.s. to $\int x^2\alpha_i(dx)$, if $i=j$, and to $\int x \alpha_i(dx)\int x \alpha_j(dx)$, if $i\neq j$. On the other hand,
	by Propositions \ref{th:convpred} and  \ref{prop:suffcond_r}, $\mathbb E[X_{n,i}\mid\mathcal G_{n-1}]\mathbb E[X_{n,j}\mid\mathcal G_{n-1}]$ converges $\mathbb P$-a.s. to $\int x\alpha_i(dx)\int x\alpha_j(dx)$, as $n\rightarrow\infty$. Thus, $\frac{1}{n}\sum_{k=1}^{n} \mathbb E[X_{k,i}\mid\mathcal G_{k-1}] \, \mathbb E[X_{k,j}\mid\mathcal G_{k-1}]$ converges $\mathbb P$-a.s. to the same limit. Hence, $\frac{1}{n}\sum_{k=1}^{n} U_{k,i} U_{k,j}$ converges in probability to $\sigma_{\alpha_i}^2$ if $j=i$ and to zero otherwise. By (\ref{eq:Z^2}), $\sum_{k=1}^{n}Z_{n,k}^2$ converges in probability to $\sum_{i,j\in J}a_ia_j\Sigma_{i,j}$.\\
\noindent $\square$
	
The next result investigates the asymptotic behavior of the deviation of the sample mean from the conditional expectation. The assumptions involve the forecast errors $U_{n,i}$, defined in (\ref{eq:Uki}), and the \textit{prediction increments}
	\begin{equation}
	\label{eq:Delta}
	\Delta E_{n,i}=\mathbb{E}\left[X_{n+1,i}\mid\mathcal{G}_n\right]- \mathbb{E}\left[X_{n,i}\mid\mathcal{G}_{n-1}\right]
	, \quad  n\geq 1,i\in I.
	\end{equation}


	\begin{thm}\label{thm:clt}
		Let $[X_{n,i}]$ be a partially $\calG$-c.i.d.
 array of real-valued and square integrable random variables and let
			\begin{equation*}
				\tilde{S}_{n,i}:=\sqrt n\left( \overline{X}_{n,i}-\mathbb{E}[X_{n+1,i}\mid\mathcal{G}_{n}]\right), \quad  n\geq 1, i\in I,
		\end{equation*}
where
$\overline{X}_{n,i}=\sum_{k=1}^nX_{k,i} / n $. 
If 
			\begin{equation}
	\label{eq:clth1}
	\underset{n}{\sup}\,\mathbb{E}[{\tilde{S}_{n,i}}^2]<\infty \quad\quad \mbox{ for every } i\in I,
	\end{equation}
and 
		\begin{equation}
	\label{eq:clth2}
	\frac{1}{n}\sum_{k=1}^{n}V_{k,i}V_{k,j}\overset{}{\rightarrow}\tilde\Sigma(e_i,e_j), \;\;\mathbb{P}\text{-a.s.}\quad\quad\quad \mbox{ for every }i,j\in I, 
	\end{equation}
where $V_{k,i}=U_{k,i}-k\Delta E_{k,i}$, with $U_{k,i}$ and $\Delta E_{k,i}$ as in (\ref{eq:Uki}) and (\ref{eq:Delta}), then $(\mathbf{\tilde{S}}_{n}^I)$ converges stably with representing measure {\rm $\Norm(0,\tilde\Sigma)$}.
						\end{thm}
\noindent	\textsc{Proof.} We can write
	\[
	\tilde{S}_{n,i}=\mathbb E[\tilde{S}_{n,i}\mid \mathcal G_n]=\sum_{k=1}^n\left( \mathbb E[\tilde{S}_{n,i}\mid \mathcal G_k ] - \mathbb E[\tilde{S}_{n,i}\mid \mathcal G_{k-1} ]  \right).
	\]
	The sequence $(\mathbb E[\tilde{S}_{n,i}\mid \mathcal G_k ] - \mathbb E[\tilde{S}_{n,i}\mid \mathcal G_{k-1} ] )_{k\geq 1}$ is a martingale difference with respect to $(\mathcal G_k)_{k\geq 1}$. Furthermore,
	\[
	\mathbb E[\tilde{S}_{n,i}\mid\mathcal G_k]-\mathbb E[\tilde{S}_{n,i}\mid\mathcal G_{k-1}]=U_{k,i}-k\Delta E_{k,i}=V_{k,i}.
	\]
Therefore,  $\tilde{S}_{n,i}=n^{-1/2} \sum_{k=1}^n V_{k,i}$. Let us fix a finite set $J\subset I$ and real numbers  $(a_i,i\in J)$, and prove that $\sum_{i\in J}a_i \tilde{S}_{n,i}$ converges stably with representing measure $\Norm(0,\sum_{i,j\in J}a_{i}a_j\tilde\Sigma_{i,j})$, where $\tilde\Sigma_{i,j}=\tilde\Sigma(e_i,e_j)$.\\
	We use Theorem \ref{lemma:clt}, with $\mathcal{F}_{n,k}=\mathcal{G}_k$ and $Z_{n,k}=n^{-1/2}\sum_{i\in J}a_iV_{k,i}$. Then $\sum_{i\in J}a_i \tilde{S}_{n,i}=\sum_{k=1}^n Z_{n,k}$.
Conditions \eqref{enum:clt_1} and \eqref{enum:clt_2} are immediate and \eqref{enum:clt_5} follows from \eqref{eq:clth2}, by noticing that
	\[
	\sum_{k=1}^{n}Z_{n,k}^2=\frac{1}{n}\sum_{k=1}^{n}\left(\sum_{i\in J}a_i V_{k,i}\right)^2 \rightarrow \sum_{i,j\in J}a_ia_j\tilde\Sigma_{i,j}\quad\mathbb{P}\text{-a.s.}.
	\]
	To prove  \eqref{enum:clt_3} and \eqref{enum:clt_4}, notice that
	$$
	Z_{nn}^2=\sum_{k=1}^n\frac{\sum_{i\in J}(a_iV_{k,i})^2}{n}-\frac{n-1}{n}\sum_{k=1}^{n-1}\frac{(\sum_{i\in J}a_iV_{k,i})^2}{n-1}\rightarrow 0\;\;\mathbb{P}\text{-a.s.}
	$$
	and
	$ Z_{n,k}=Z_{k,k}{\sqrt{k/n}}	$.
	Thus, $\underset{1\leqslant k \leqslant n} {\max}Z_{n,k}^2 {\rightarrow}0$, $\mathbb P$-a.s. and
	$$
	\underset{n}{\sup}\mathbb{E}[\underset{1\leqslant k\leqslant n}{\max}Z_{n,k}^2]\leqslant\underset{n}{\sup} \sum_{k=1}^n\mathbb{E}\left(Z_{n,k}^2\right)\leqslant \sum_{i,j\in J}|a_ia_j|\;\underset{n}{\sup}\mathbb{E}|\tilde{S}_{n,i}\tilde{S}_{n,j}|<\infty.
	$$
Thus, by Theorem \ref{lemma:clt}, the thesis follows.	\\
\noindent $\square$

\medskip

 \begin{exa}\rm {\bf [Interacting randomly reinforced processes (Ctd)]}.
 	As an illustration of Theorems \ref{thm:clt1} and \ref{thm:clt}, let us	consider again the interacting randomly reinforced process of Example  \ref{ex:interactingRRprocesses}. Assume that the $X_{n,i}$ are generated independently as in (\ref{eq:recursive}).
 	Interaction is given by common weights $W_n$, where $W_n \sim p_W$ independently from the past and from the concomitant values $X_{n,i}$. Let us assume that the common distribution $p_W$ has finite second moment, and that the $X_{n,i}$ are square-integrable, $\mathbb R$-valued random variables. Hence, $(\bX_n^I)$
 	is partially $\mathcal G$-c.i.d., where  $\calG_n=\sigma(\bX_{1:n}^I, \mathbf W_{1:n}^I)$. 
 By Proposition \ref{prop:suffcond_cid}, $(X_{n,i})_{n\geq 1}$ and $(X_{n,i}^2)_{n\geq 1}$ are $\mathcal G$-predictive uniformly integrable and dominated in $L^1$. Hence, $\mathbb E[X_{n+1,i}\mid\mathcal G_n]\rightarrow \int x \alpha_i(dx)<\infty$ and $\mathbb E[X_{n+1,i}^2\mid\mathcal G_n]\rightarrow \int x^2 \alpha_i(dx)<\infty$, $ \mathbb P$-a.s., 
 	where $(\alpha_i:i\in I)$ are the directing random measures of the array $[X_{n,i}]$. 
 	Let us denote
 	$\mu_{\alpha_i}=\int xd\alpha_i$ and $\sigma_{\alpha_i}^2=\int x^2\alpha_i(dx)-\mu_{\alpha_i}^2$;
furthermore, let
$\Phi(\cdot;0, \sigma^2)$ be the distribution function of the Gaussian law with parameters $0$ and $\sigma^2$, with the proviso that
$\Phi(\cdot;0,0)$ represents the distribution degenerate at $0$.\\
By Theorem \ref{thm:clt1}  and the properties of stable convergence, for every finite $J\subset I$
 and every $s_j\neq 0$, $j\in J$,
	\[
	\mathbb P\left[\bigcap_{i\in J}\left(\sum_{k=1}^n(X_{k,i}-\mathbb E[X_{k,i}\mid \mathcal G_{k-1}])\leq {\sqrt{n}}s_i\right)\right]\rightarrow
	\int \prod_{i\in J}\Phi(s_i;0,\sigma^2_{\alpha_i}(\omega))\mathbb P(d\omega),\quad\quad \mbox{ as }n\rightarrow\infty.
	\]
Moving to Theorem \ref{thm:clt},
notice that \[
U_{n,i}=\left(X_{n,i}-\mu_{n,i}\right),\quad \Delta E_{n,i}=\frac{\left(X_{n,i}-\mu_{n,i}\right)W_n}{w_0+\sum_{k=1}^nW_k},\quad
V_{n,i}=\left(X_{n,i}-\mu_{n,i}\right)\left(1-\frac{W_n}{w_0/n+\sum_{k=1}^nW_k/n}\right),
 	\]
 	where $\mu_{n,i}=\mathbb E[X_{n,i}\mid\mathcal G_{n-1}]$ . Thus, 	
 	\[
 	\begin{aligned}
 	\sup_n\mathbb E[\tilde{S}_{n,i}^2]\leq \sup_n\mathbb E[V_{n,i}^2] &\leq
 	\sup_n\mathbb E[X_{n,i}^2] \left(1+\sup_n \mathbb E\left[\left(\frac{1}{n}+\frac{1}{n}\sum_{k=1}^{n-1}\frac{W_k}{W_n}\right)^{-2}\right]\right)\\&\leq
 	\sup_n\mathbb E[X_{n,i}^2]
 	\left(1+\sup_n \mathbb E\left[ \frac{1}{n}+\frac{1}{n}\sum_{k=1}^{n-1}\frac{W_n^2}{W_k^2}\right] \right)
 	\end{aligned}\]
 	where $1/0:=\infty$ and the last inequality comes from the convexity of $f(x):=1/x^{2}$.
 	The above inequalities show that a sufficient condition for (\ref{eq:clth1}) to hold is $\mathbb E[1/W_1^2]<\infty$.
We restrict our attention to this case.
 	 	 	To prove (\ref{eq:clth2}), notice that
\[ 	\frac{1}{n}\sum_{k=1}^nV_{k,i}V_{k,j}
 	 =\frac{1}{n}\sum_{k=1}^n\left(X_{k,i}-\mu_{k,i}\right)\left(X_{k,j}-\mu_{k,j}\right)\left(1-\frac{W_k}{w_0/k+\sum_{h=1}^k{W}_h/k}\right)^2.
 	\]	
Moreover, $[W_n,\bX_{n}^I]_{n\geq 1}$ is a partially $\calG$-c.i.d. array, with directing measures  $(p_W,(\alpha_i:i\in I))$ and, proceeding as in Example \ref{exa:1}, it can be proved that the sequence of products $(X_{n,i} X_{n,j})_{n\geq 1}$ is predictive uniformly integrable and dominated in $L^1$, for every $i,j\in I$. Since $\mathbb E[X_{n,i}X_{n,j}W_n^2]=\mathbb E[X_{n,i}X_{n,j}]\mathbb E[W_n^2]$ and
 	$\mathbb E[X_{n,i} X_{n,j} W_n^2\mid \mathcal G_{n-1}]=
 	\mathbb E[X_{n,i} X_{n,j} \mid \calG_{n-1}] \mathbb E[W_n^2]$, the sequence of products $(X_{n,i} X_{n,j} W_n^2)_{n\geq 1}$ is also predictive uniformly integrable and dominated in $L^1$. The same reasoning applies to $(X_{n,i} W_n^2)_{n\geq 1}$. Thus, by Theorem \ref{th:slln}, as $n\rightarrow\infty$,
 	\begin{equation}
 	\label{eq:exa_centrlth_conv}
 	\begin{aligned}
 	 	&\frac{1}{n}\sum_{k=1}^nX_{k,i}X_{k,j}\stackrel{a.s.}{\rightarrow}\left\{ \begin{array}{ll}\mu_{\alpha_i}\mu_{\alpha_j} & j\neq i\\ (\sigma_{\alpha_i}^2+\mu_{\alpha_i}^2)&j=i\end{array}\right. \\ &\frac{1}{n}\sum_{k=1}^n X_{k,i}X_{k,j}W_k^2\stackrel{a.s.}{\rightarrow} \left\{\begin{array}{ll} \mu_{\alpha_i}\mu_{\alpha_j}\mathbb E[W_1^2] &j\neq i\\ (\sigma_{\alpha_i}^2+\mu_{\alpha_i}^2)\mathbb E[W_1^2] &j=i\end{array}\right.\\
 	&\frac{1}{n}\sum_{k=1}^n X_{k,i}W_k^2\stackrel{a.s.}{\rightarrow} \mu_{\alpha_i}\mathbb E[W_1^2],\end{aligned}
 	\end {equation}
 	where $\mu_{\alpha_i}=\int x \alpha_i(dx)$ and $\sigma_{\alpha_i}^2=\int x^2\alpha_i(dx)-\mu_{\alpha_i}^2$.
 	Noticing that $\mu_{n,i}\stackrel{a.s.}{\rightarrow}\mu_{\alpha_i}$ and $1/(w_0/n+\sum_{k=1}^n W_k/n)\stackrel{a.s.}{\rightarrow} 1/\mathbb E[W_1]$,  we obtain
 	\[
 	 	\label{eq:exa_centrlth}
 	\frac{1}{n}\sum_{k=1}^nV_{k,i}V_{k,j}
 	 	\stackrel{a.s.}{\rightarrow}\left\{
 	\begin{array}{ll}
 	\sigma_{\alpha_i}^2{\mathbb V(W_1)}/{\mathbb E[W_1]^2} & j=i\\
 	0 & j\neq i,
 	\end{array}\right.
 	 	\]
 	as $n\rightarrow\infty$. In the above equation, $\mathbb E[W_1]>0$, since we are assuming $\mathbb E[1/W_1^2]<\infty$. \\
 	 	Therefore, the assumptions of Theorem \ref{thm:clt} are satisfied, and, for every finite $J\subset I$ and every $s_i\neq 0$ ($i\in J)$,
 	\[
 	\mathbb P\left[\bigcap_{i\in J}\left(\overline{X}_{n,i}-\mathbb E[X_{n+1}\mid\mathcal G_n]\leq \frac{s_i}{\sqrt{n}}\right)\right]\rightarrow
 	\int \prod_{i\in J}\Phi\left(s_i;0,\sigma_{\alpha_i}^2(\omega)\frac{\mathbb V[W_1]}{\mathbb E[W_1]}\right)\mathbb P(d\omega).
 	\]
  Notice that, if $\mathbb V[W_1]=0$   (i.e. if the $(X_{n,i})_{n\geq 1}$ are independent sequences of exchangeable random variables), then $\sqrt n(\overline X_{n,i}-\mathbb E[X_{n+1}\mid\mathcal G_{n}])$ converges in probability to zero. \\
  \noindent $\square$
  	
\end{exa}

In the above example, $S_n$ and $\tilde{S}_n$ have a similar behaviour. In particular, the limiting covariance matrices are both diagonal. This is not always the case, as the next example shows.

\begin{exa}
\rm
Let us consider the partially c.i.d. array $[X_{n,i}]_{n\geq 1,i=1,2}$ defined in Example \ref{ex:uniform}.
Since the $X_{n,i}$'s are bounded, they are square integrable. Hence,  $\mu_{\alpha_i}=\int x\alpha_i(dx)<\infty$ and  $\sigma_{\alpha_i}^2=\int x^2\alpha_i(dx)-\mu_{\alpha_i}^2<\infty$. By Theorem \ref{thm:clt1}, for every $s_1,s_2\neq 0$,
$$\mathbb P\left[\bigcap_{i=1}^2\left[\frac{\sum_{k=1}^n(X_{k,i}-\mathbb E[X_{k,i}\mid\mathcal G_{k-1}])}{\sqrt n}\leq s_i\right]\right]\rightarrow
\int \prod_{i=1}^2 \Phi(s_i;0,\sigma^2_{\alpha_i}(\omega))\mathbb P(d\omega).$$
Let us now turn to Theorem \ref{thm:clt}.
A trivial calculation shows that, for every $k$ and $i$,
\[
V_{k,i}
=(1-k\beta_kX_{k,j})\left(X_{k,i}-\mathbb E[X_{k+1,i}\mid\mathcal G_k]\right).
\]
The above equation suggests that $\sqrt n(\overline X_{n,i}-\mathbb E[X_{n+1,i}\mid\mathcal G_n])$ may fail to converge, for specific settings of $(\beta_n)$ (e.g. if $\beta_n=1$ for every $n$). Here, we consider $\beta_n=2/(n+1)$.
Condition (\ref{eq:clth1}) holds, since $(V_{n,i})_{n\geq 1}$ is bounded and
\[
\sup_n\mathbb E[\tilde{S}_{n,i}^2]\leq \sup_n\mathbb E[V_{n,i}^2]<\infty.
\]
Furthermore, $\mathbb P$-a.s.
\[
\frac{1}{n}\sum_{k=1}^nV_{k,i}V_{k,j}
=\frac{1}{n}\sum_{k=1}^n\left(X_{k,i}-\mu_{k,i}\right)\left(X_{k,j}-\mu_{k,j}\right)\left(1-\frac{2k}{k+1}X_{k,i}\right)\left(1-\frac{2k}{k+1}X_{k,j}\right)
\rightarrow\tilde\Sigma_{i,j},\,\,\mathbb P\mbox{-a.s.},
\]
where \begin{equation*}
\tilde\Sigma=4
\left[
\begin{array}{cc}
\int (x-\mu_{\alpha_1})^2(x-1/2)^2\alpha_1(dx)& \sigma_{\alpha_1}^2\sigma_{\alpha_2}^2\\
\sigma_{\alpha_1}^2\sigma_{\alpha_2}^2&\int (x-\mu_{\alpha_2})^2(x-1/2)^2\alpha_2(dx)
\end{array}
\right].
\end{equation*}
According to Theorem \ref{thm:clt},
$$\mathbb P\left[\bigcap_{i=1}^2\left(\overline X_{n,i}-\mathbb E[X_{n+1,i}\mid\mathcal G_n]\leq \frac{s_i}{\sqrt n}\right)\right]\rightarrow\int \Phi(s_1,s_2;0,\tilde\Sigma(\omega))\mathbb P(d\omega).$$
\noindent $\square$
\end{exa}

\noindent
{\bf Acknowledgments}. We are grateful to Pietro Rigo for helpful suggestions. We sincerely thank the anonymous referees for their very accurate reviews. The authors have been partially supported by grants from Bocconi University and PRIN grant 2015SNS29B.

\bibliographystyle{plainnat}   
\bibliography{pcid}

\section{Appendix.}  \label{sec:appendix}

Here we provide some complements for Section \ref{sec:partialExch}.


\begin{prop}
Consider the Gaussian predictive system described  in Example \ref{ex:gauss}. 
Then:
\begin{itemize}
	\item[(i)] the process $(\bX_n^I)$ is $\mathcal G$-partially c.i.d. with directing measures {\rm $\alpha_i = \Norm(\mu_i, \gamma \sigma^2_{1,i})$}, where $\gamma=\prod_{k=1}^\infty (1-\lambda_k^2)$ and the  random means $\mu_i$ are conditionally independent, given $\gamma$, with 
{\rm $\mu_i \mid \gamma \stackrel{ind}{\sim}\Norm(\mu_{1,i},(1-\gamma)\sigma_{1,i}^2)$, $i=1, \ldots, K$}.
	\item[(ii)] The random probability measures $\alpha_1,\dots,\alpha_K$ are not stochastically independent  when the observations arrive at a Poisson process rate.
	\end{itemize}
\end{prop}

\noindent {\sc Proof} \\
\noindent (i) For showing that $X_{n+2,i}\mid\mathcal G_n^i\stackrel{d}{=}X_{n+1,i}\mid\mathcal G_n^i$, let us work with the characteristic function, and notice that  
\[
\begin{aligned}
\mathbb E[\exp(isX_{n+2,i})\mid\mathcal G_n^i] 
&= \exp(is((1-\lambda_{n+1})\mu_{n+1,i}-s^2(1-\lambda_{n+1}^2)\sigma_{n+1,i}^2/2)\mathbb E[\exp(is\lambda_{n+1}X_{n+1,i}) \mid\mathcal G_{n}^i]\\
&=\exp(is(\mu_{n+1,i}-s^2\sigma_{n+1,i}^2/2))
\end{aligned}
\]
Therefore, $(\bX_n^I)$ is partially $\calG$-c.i.d. By Proposition \ref{prop:strong}, the marginal directing measure $\alpha_i$ ìs the limit of the predictive distribution $P_{n,i}$, which is a $\Norm(\mu_{n,i}, \sigma^2_{n,i})$ as in (\ref{eq:gaussPred}). For each $i$, the 
sequence $(\mu_{n,i})$ is a uniformly integrable martingale. Hence, $\mu_{n,i}$ converges to an integrable random variable $\mu_i$, $\mathbb P$-a.s. As for the variances, the sequence $(\sigma^2_{n,i})$ is decreasing an non negative. Hence, it converges $\mathbb P$-a.s. to a non negative random variable $\sigma^2_i$. By the recursive equations (\ref{eq:gaussrecurs}), one obtains 
$\sigma^2_i=\gamma \sigma^2_{1,i}$, where  $\gamma=\prod_{k=1}^\infty (1-\lambda^2_k)$. 
Thus, $P_{n,i} \rightarrow \Norm(\mu_i, \sigma_i^2)$. 

To find the probability distribution of the random vector $\mu=(\mu_1, \ldots, \mu_K)$, notice that, 
conditionally on $\lambda:=(\lambda_n)$, $(X_{1:n})$ has a Gaussian distribution.  Therefore, given $\lambda$, the vector $\mu_{n}=(\mu_{n,1}, \ldots, \mu_{n,K})$ is Gaussian, with expectation $\mu_{1}=(\mu_{1,1},\dots,\mu_{1,K})$. Furthermore, since $\mu_{n,1},\dots,\mu_{n,K}$ are conditionally independent, given $\lambda$, then $\mu_{1},\dots,\mu_{K}$ are also conditionally independent given $\lambda$. To compute the conditional variances,  notice that
\[
\begin{aligned}
\mathbb E[\mu_{n,i}^2\mid \lambda]=\mathbb E\left[ \left((1-\lambda_{n-1})\mu_{n-1,i}+\lambda_{n-1} X_{n-1,i}   \right)^2\mid\lambda\right]
=(1-\lambda_{n-1}^2)\mathbb E[\mu_{n-1,i}^2\mid\lambda]+\lambda_{n-1}^2(\mu_{1,i}^2+\sigma_{1,i}^2)
\end{aligned}
\]
Given the initial value $\mathbb E[\mu^2_{1,i}\mid\lambda]=\mu_{1,i}^2$, the unique solution of the above system is
\[
\mathbb E[\mu_{n,i}^2\mid \lambda]=\mu_{1,i}^2+\left[1-\prod_{k=1}^{n-1}(1-\lambda_k^2)\right]\sigma_{1,i}^2.\]
Thus,
\[
\mathbb V[\mu_{n,i}\mid\lambda]=\left[1-\prod_{k=1}^{n-1}(1-\lambda_k^2)\right]\sigma_{1,i}^2.
\]
Since $\mu_i$ is the limit of $\mu_{n,i}$,
\[
\mu_i\mid\lambda\sim \Norm(\mu_{1,i},(1-\gamma)\sigma_{1,i}^2),\]
with $\gamma=\prod_{k=1}^\infty (1-\lambda_k^2)$. Since the conditional distribution of $(\mu_1,\dots,\mu_K)$, given $\lambda$, is a function of $\gamma$, then (i) holds.

\medskip
\noindent (ii)
The random probability measures $\alpha_1,\dots,\alpha_K$ are conditionally independent, given $\gamma$, but are not stochastically independent, unless $\gamma$ has a degenerate distribution. For example,
$Cov(\mu_i^2,\mu_j^2)=\sigma_{1,i}^2\sigma_{1,j}^2\mathbb V(\gamma)$.\\
When the observations arrive at a Poisson process rate, $\lambda_1,\lambda_2,\dots$ are independent random variables, and $\lambda_n\sim\beta(1,n)$. Thus
\[
\mathbb E[\gamma]=\lim_{n\rightarrow\infty}\prod_{k=1}^n\left(
1-\mathbb E[\lambda_k^2]\right)=\lim_{n\rightarrow\infty}\prod_{k=1}^n\left(
1-\frac{2}{(k+1)(k+2)}
\right)=1/3
\]
while
\[\begin{aligned}
\mathbb V[\gamma]&=\lim_{n\rightarrow\infty}\left(\prod_{k=1}^n
\mathbb E[\lambda_k^4-2\lambda_k^2+1]-\prod_{k=1}^n\left(
1-\mathbb E[\lambda_k^2]\right)\right)\\
&=\lim_{n\rightarrow\infty}\left(\prod_{k=1}^n\left(1-\frac{4}{(k+1)(k+2)}+\frac{24}{(k+1)(k+2)(k+3)(k+4)}\right)
-\prod_{k=2}^n\left(1-\frac{2}{(k+1)(k+2)}\right)^2\right)\\
&\geq \frac{4}{45}\lim_{n\rightarrow\infty}\prod_{k=2}^n\left(1-\frac{2}{(k+1)(k+2)}\right)^2>0.
\end{aligned}\]
Hence $\gamma$ is not constant, therefore $\alpha_i$ is not degenerate and $\alpha_1,\dots,\alpha_K$ are not stochastically independent. \\
\noindent $\square$

\end{document}